\documentclass[english,10pt]{article}

\usepackage{xcolor}
\usepackage{amsmath}
\usepackage{graphics}
\usepackage{epsfig}
\usepackage{vmargin}
\usepackage{eqnarray}
\usepackage{amssymb}
\usepackage[latin1]{inputenc}
\usepackage{babel,indentfirst,amsfonts,array}

\newtheorem{thm}{Theorem}[section]
\newtheorem{lemme}[thm]{Lemma}
\newtheorem{prop}[thm]{Proposition}
\newtheorem{cor}[thm]{Corollary}

\newtheorem{defi}[thm]{Definition}

\def\P{\mathbb P}
\def\E{\mathbb E}

\def\R{\mathbb R}
\def\C{\mathbb C}

\def\Z{\mathbb Z}
\def\N{\mathbb N}
\def\({\left(}
\def\){\right)}
\def\[{\left[}
\def\]{\right]}
\def\fin{\hfill\square}

\def\fin{\hfill $\square$}

\let\oldsqrt\sqrt
\def\sqrt{\mathpalette\DHLhksqrt}
\def\DHLhksqrt#1#2{%
\setbox0=\hbox{$#1\oldsqrt{#2\,}$}\dimen0=\ht0
\advance\dimen0-0.2\ht0
\setbox2=\hbox{\vrule height\ht0 depth -\dimen0}%
{\box0\lower0.4pt\box2}}
\title{\textsc{Random walk in a stratified medium}
\author{Julien Br\'emont}
\date{Universit\'e Paris-Est Cr\'eteil,~mai 2016}
}

\begin{document}

\maketitle

\setcounter{page}{1}

\begin{abstract}
We give a recurrence criterion for a Markov chain in $\Z^{d+1}$ in a medium stratified by parallel affine hyperplanes. The asymptotics of the random walk is governed by some notion of directional flux variance, describing the dispersive power of some associated average flow. The result admits a geometrical interpretation, surprisingly intrinsically non-Euclidean. Some applications and open questions are discussed. 
\end{abstract}

\footnote{
\begin{tabular}{l}\textit{AMS $2000$ subject classifications~: 60J10, 60K20} \\
\textit{Key words and phrases~: Markov chain, recurrence criterion, continued fraction, anisotropic pseudosphere.} 
\end{tabular}}

\section{Introduction}
We study the recurrence properties of an inhomogeneous Markov chain $(S_n)_{n\geq0}$ in $\Z^{d}\times\Z$, where $d\geq1$. Starting the random walk at $0$, let $S_n=(S_n^1,S_n^2)\in\Z^d\times\Z$. We call ``vertical'' the quantities relative to the second coordinate. The environment is invariant under $\Z^d$-translations, i.e. the collection of transitions laws is stratified with respect to the affine hyperplanes $(\Z^d\times\{n\})_{n\in \Z}$. We make no hypothesis on the relative dependence between transitions laws in distinct hyperplanes. 

\smallskip
\noindent
A planar random walk of this type was proposed by Campanino and Petritis \cite{cp} in 2003, as a simplified probabilistic version of PDE transport models in stratified porous medium considered by Matheron and de Marsily \cite{mama}. Following this line of research, we focus on a more general case in $\Z^{d+1}$. For the sequel we fix Euclidean Norms and denote scalar product by a dot. 

\medskip
Let us state the model. For each vertical $n\in\Z$, let reals $p_n,q_n,r_n$ with $p_n+q_n+r_n=1$ and a probability measure $\mu_n$ with support in $\Z^d$. We suppose that for some $\delta>0$ and all $n\in\Z$~:

\smallskip
1) $\min\{p_n,q_n,r_n\}\geq\delta$,

\smallskip
2) $\sum_{k\in\Z^d}\|k\|^{\max(d,3)}\mu_n(k)\leq1/\delta$,

\smallskip
3) the eigenvalues of the real symmetric matrix $\sum_{k\in\Z^d}kk^T~\mu_n(k)$ are $\geq\delta$. Equivalently~:

$$\sum_{k\in\Z^d}(t.k)^2\mu_n(k)\geq \delta \|t\|^2,~t\in\R^d.$$

\noindent
Notice that the last condition implies that the subgroup of $(\Z^d,+)$ generated by $\mbox{supp}(\mu_n)$ is $d$-dimensional. The transition laws are then defined, for all $(m,n)\in\Z^d\times\Z$ and $ k\in\Z^d$, by~:

$$\P_{(m,n),(m,n+1)}=p_n,~\P_{(m,n),(m,n-1)}=q_n,~\P_{(m,n),(m+k,n)}=r_n\mu_n(k).$$

\begin{center}
\resizebox{0.42\linewidth}{!}{\begin{picture}(0,0)%
\includegraphics{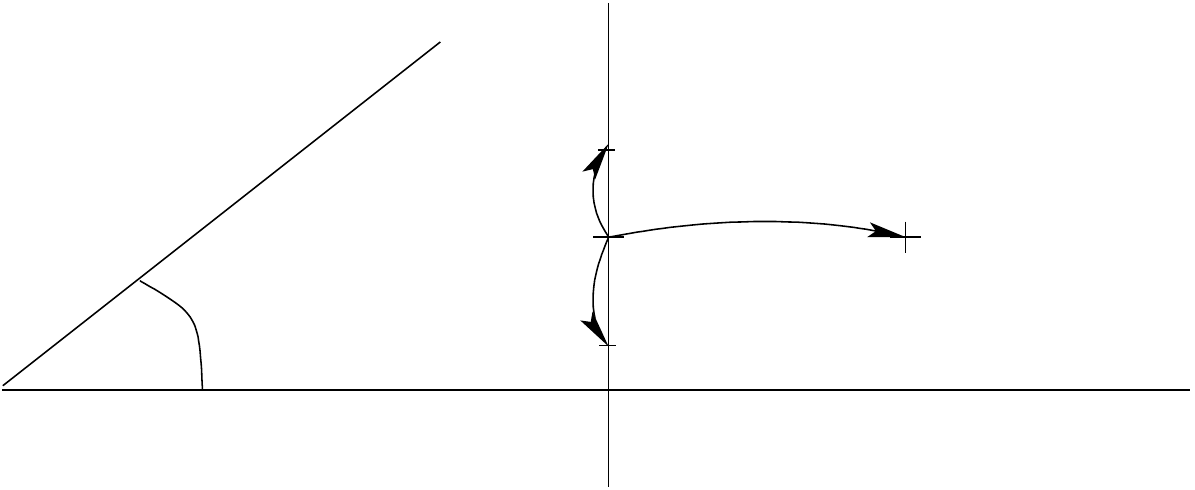}%
\end{picture}%
\setlength{\unitlength}{3947sp}%
\begingroup\makeatletter\ifx\SetFigFont\undefined%
\gdef\SetFigFont#1#2#3#4#5{%
  \reset@font\fontsize{#1}{#2pt}%
  \fontfamily{#3}\fontseries{#4}\fontshape{#5}%
  \selectfont}%
\fi\endgroup%
\begin{picture}(5724,2349)(4279,-4348)
\put(7276,-2161){\makebox(0,0)[lb]{\smash{{\SetFigFont{12}{14.4}{\rmdefault}{\mddefault}{\updefault}{\color[rgb]{0,0,0}$\Z$}%
}}}}
\put(7501,-2911){\makebox(0,0)[lb]{\smash{{\SetFigFont{12}{14.4}{\rmdefault}{\mddefault}{\updefault}{\color[rgb]{0,0,0}$r_n\mu_n(k)$}%
}}}}
\put(6826,-2911){\makebox(0,0)[b]{\smash{{\SetFigFont{12}{14.4}{\rmdefault}{\mddefault}{\updefault}{\color[rgb]{0,0,0}          $p_n$}%
}}}}
\put(8701,-3286){\makebox(0,0)[lb]{\smash{{\SetFigFont{12}{14.4}{\rmdefault}{\mddefault}{\updefault}{\color[rgb]{0,0,0}$(m+k,n)$}%
}}}}
\put(6751,-3436){\makebox(0,0)[lb]{\smash{{\SetFigFont{12}{14.4}{\rmdefault}{\mddefault}{\updefault}{\color[rgb]{0,0,0}$q_n$}%
}}}}
\put(6533,-3196){\makebox(0,0)[lb]{\smash{{\SetFigFont{12}{14.4}{\rmdefault}{\mddefault}{\updefault}{\color[rgb]{0,0,0}$(m,n)$}%
}}}}
\put(4666,-3811){\makebox(0,0)[lb]{\smash{{\SetFigFont{12}{14.4}{\rmdefault}{\mddefault}{\updefault}{\color[rgb]{0,0,0}$\mathbb{Z}^d$}%
}}}}
\end{picture}%
}
\end{center}

\medskip
The model of Campanino-Petritis corresponds to taking $d=1$, with $p_n=q_n=p\in(0,1)$ and $\mu_n=\delta_{\varepsilon_n}$, fixing some sequence $(\varepsilon_n)_{n\in\Z}$ of $\pm1$. Campanino and Petritis \cite{cp} for instance show recurrence when $\varepsilon_n=(-1)^n$ and transience for $\varepsilon_n=1_{n\geq 0}-1_{n<0}$ or when the $(\varepsilon_n)$ are typical realizations of i.i.d. random variables with law $(\delta_1+\delta_{-1})/2$. In some neighbourhood of this setting, several variations, extensions and second order questions were subsequently considered by various authors; see the introduction of \cite{jb1}. In \cite{jb1}, a recurrence criterion was given for the model introduced above when $d=1$, assuming the local vertical symmetries $p_n=q_n$, $n\in\Z$. In this family of random walks, planar simple random walk, hardly recurrent, is the most recurrent one. This explains the prevalence of transience results on the Campanino-Petritis model. Mention that for the latter, a growth condition larger than $\log n$ on $\varepsilon_1+\cdots+\varepsilon_n$ is sufficient to ensure transience.

\medskip
Pushing to some natural limit the method used in \cite{jb1}, we establish in this article a recurrence criterion for the model described above. This furnishes a large class of recurrent random walks in $\Z^2$ and $\Z^3$. The mechanism governing the asymptotic behaviour of the random walk reveals some familiarity with classical Electromagnetism, involving notions such as flux variations. The latter represent the dispersive properties of some average flow associated with the random walk. Variations are measured in a probabilistic sense, via some empirical variances. We also provide a geometrical interpretation of the recurrence criterion. Surprisingly it involves hyperbolic geometry, stereographic projections and some kind of anisotropic pseudosphere.

\section{Statement of the result}
\subsection{Notations and result}
\begin{defi}

 $ $
\begin{itemize}
\item For $n\in\Z$, let $m_n=\sum_{k\in\Z^d}k\mu_n(k)$ be the expectation of $\mu_n$. 

\item
For $n\in\Z$, let $p'_n=p_n/(p_n+q_n)$, $q'_n=q_n/(p_n+q_n)$.

\item For $n\in\Z$, let $a_n=q'_n/p'_n=q_n/p_n$ and $b_n=1/p'_n=1+a_n$. 

\item Set ~:

$$\rho_n=\left\{{\begin{array}{cc}
a_1\cdots a_{n},&n\geq1,\\
1,&n=0,\\
(1/a_{n+1})\cdots(1/a_{-1})(1/a_0)&~n\leq-1.\end{array}
}\right.$$

\item For $n\geq 0$, let~:

$$v_+(n)=\sum_{0\leq k\leq n}\rho_k\mbox{ and }v_-(n)=a_0\sum_{-n-1\leq k\leq-1}\rho_k,$$

\noindent
as well as~:

$$w_+(n)=\sum_{0\leq k\leq n}(1/\rho_k)\mbox{ and }w_-(n)=(1/a_0)\sum_{-n-1\leq k\leq-1}(1/\rho_k).$$

\end{itemize}

\end{defi}

We denote by $\theta$ the ``left shift" on indices. Given $f=f((q_i/p_i)_{i\in \Z})$, set $\theta f=f((q_{i+1}/p_{i+1})_{i\in \Z})$. In particular the cocycle relation for $(\rho_n)$ reads as~:

$$\forall (n,k)\in\Z^2,~\rho_{n+k}=\rho_n\theta^n\rho_k.$$

\medskip
We next need a definition of inverse function for non-decreasing functions defined on the set of non-negative integers $\N=\{0,1,\cdots\}$ and having values in $\R_+\cup\{+\infty\}$.

\begin{defi}
$ $

\noindent
Let $f:\N\rightarrow\R_+\cup\{+\infty\}$, non-decreasing. For $x\in\R$, let $f^{-1}(x)=\sup\{n\in\N~|~f(n)\leq x\}$, with $\sup\{\N\}=+\infty$ and $\sup\{\o\}=0$. 
\end{defi}

We next turn to notions related to directional fluxes and their variations.

\begin{defi}
$ $

\begin{itemize}

\item Let $S_+^{d-1}=\{x\in\R^d~|~\|x\|=1,~x_1\geq0\}$ be a half unit Euclidean sphere of $\R^d$. 

\item For $u\in S_+^{d-1}$ and $k\leq l$ in $\Z$, introduce~:

$$R_k^l(u)=\sum_{s=k}^l\frac{r_s}{p_s}\frac{\rho_l}{\rho_s}m_s.u\mbox{ and }T_k^l(u)=\frac{\rho_{k-1}}{\rho_{l}}(R_k^l(u))^2=\rho_{k-1}\rho_l\({\sum_{s=k}^l\frac{r_s}{p_s\rho_s}m_s.u}\)^2.$$

\item For $m\geq 0$, $n\geq0$, let $\psi(-m,n)$ be the positive (maybe $+\infty$) quantity such that~: 

$$\psi^2(-m,n)=nw_+\circ v_+^{-1}(n)+mw_-\circ v_-^{-1}(m).$$

\noindent
We also set $\psi(n)=\psi(-n,n)$, $\psi_+(n)=\psi(0,n)$, $\psi_-(n)=\psi(-n,0)$, for $n\geq0$.

\item For $u\in S_+^{d-1}$, $m\geq0$, $n\geq0$, let $\varphi_u(-m,n)$ be the positive (maybe $+\infty$) quantity such that~:

$$\varphi^2_{u}(-m,n)=\psi^2(-m,n)+\sum_{-v_-^{-1}(m)\leq k\leq l\leq v_+^{-1}(n)}T_k^l(u).$$

\noindent
Set for $n\geq 0$, $\varphi_u(n)=\varphi_u(-n,n)$ and $\varphi_{u,+}(n)=\psi^2(-n,n)+\sum_{-v_-^{-1}(n)\leq k\leq l\leq v_+^{-1}(n),kl>0}T_k^l(u).$ Introduce also~:

$$\varphi_{u,++}(n)=\psi^2(0,n)+\sum_{1\leq k\leq l\leq v_+^{-1}(n)}T_k^l(u)\mbox{ and }\varphi_{u,+-}(n)=\psi^2(-n,0)+\sum_{-v_-^{-1}(n)\leq k\leq l\leq-1}T_k^l(u).$$

\end{itemize}
\end{defi}

The aim of the article is to prove the following result.

\begin{thm}

\label{princip}

$ $

\noindent
The random walk is recurrent if and only if~:

$$\sum_{n\geq1}n^{-d-1}\int_{S_+^{d-1}}\frac{(\varphi_u^{-1}(n))^2}{\varphi_{u,+}^{-1}(n)}~du=+\infty.$$

\end{thm}

\subsection{Geometrical interpretation; corollaries}
Let us detail a geometrical interpretation of the above result. What comes out of the computations is the integral~:

$$\int_{u\in S_+^{d-1},0<t<1}\frac{(\varphi_u^{-1}(1/t))^2}{\varphi_{u,+}^{-1}(1/t)}~t^{d-1}dudt.$$

\noindent
It will be explained later why this quantity has the same order as the one appearing in the statement of the theorem. The term $(\varphi_u^{-1}(1/t))^2/\varphi_{u,+}^{-1}(1/t)$ essentially comes from a stereographic projection. We draw below a picture when $d=2$ (hence in $\R^3$) showing that the previous integral is the volume of some anisotropic version of Beltrami's pseudosphere. The classical pseudosphere is a model in $\R^3$ of a part of the hyperbolic plane (the whole hyperbolic plane cannot be represented in $\R^3$; theorem of Hilbert, 1901). Here is a way of visualizing this integral ($d=2$)~:

\smallskip
- Draw the vertical line passing at $0$, directed by $e_3$, the third vector of the canonical basis of $\R^3$. Fix $u\in S_+^1$. Let $P_0=u^{\bot}$ be the vectorial plane orthogonal to $u$. For $0<t<1$, let $P_t$ be the affine plane parallel to $P_0$ and passing through $tu$.

- We parametrize points on the left half of $P_0$ in polar coordinates $\rho e^{i\alpha}$, with $0\leq \alpha\leq\pi$ and $\rho>0$, as shown, starting from the Northern part of the vertical axis and turning counterclockwise.

- At each $\rho e^{i\alpha}$ we plug in direction $u$ (therefore orthogonally to the plane $P_0$) the length $1/\varphi_u(-\rho\sin(\alpha/2),\rho\cos(\alpha/2))$. When $0\leq \alpha\leq\pi$ and $\rho>0$, the $[-\rho\sin(\alpha/2),\rho\cos(\alpha/2)]$ describe all the intervals in the vertical direction containing the point $0$. We hence obtain a surface above the left half of $P_0$ in direction $u$, parametrized by $\rho>0$ and $0\leq \alpha\leq\pi$.

- For $0<t<1$, the plane $P_t$ cuts this surface along some dashed line shown on the picture. The point on this line and in the horizontal plane has coordinates $(t,-\sqrt{2}\varphi_u^{-1}(1/t),0)$ in the basis $(u,u',e_3)$, where $u'=e_3\wedge u$. The point on the level line lying on the vertical line passing through $(t,0,0)$ has components $(t,0,\varphi_{u,++}^{-1}(1/t))$.

- A little of geometry, related to some kind of stereographic projection, shows how to obtain $z_t$ equal to $(\varphi_u^{-1}(1/t))^2/\varphi_{u,++}^{-1}(1/t)$ up to multiplicative constants on the picture, at a point of coordinates $(t,0,z_t)$, still with respect to $(u,u',e_3)$. We use that $\varphi_u^{-1}(x)\leq\varphi_{u,+}^{-1}(x)\leq \varphi_{u,+\pm}^{-1}(x)$, giving that the orthogonal triangle in $P_t$ with vertices $(t,-\sqrt{2}\varphi_u^{-1}(1/t),0)$, $(t,0,0)$ and $(t,0,\varphi_{u,++}^{-1}(1/t))$ has a vertical edge larger or equal to the horizontal side.

- In the plane generated by $u$ and $e_3$ we have a function $t\longmapsto z_t$, $0<t<1$. The integral of this function is the top of the hatched area and has order~: 

$$\int_{0<t<1}\frac{(\varphi_u^{-1}(1/t))^2}{\varphi_{u,++}^{-1}(1/t)}~dt.$$

\noindent
Do next the same work on the Southern side and obtain a similar area (not equal to previous one in general) corresponding to~:

$$\int_{0<t<1}\frac{(\varphi_u^{-1}(1/t))^2}{\varphi_{u,+-}^{-1}(1/t)}~dt.$$

\noindent
When summing the two last integrals, one globally obtains the full hatched area, in the plane generated by $e_3$ and $u$ and this has order~:
 
 $$\int_{0<t<1}\frac{(\varphi_u^{-1}(1/t))^2}{\varphi_{u,+}^{-1}(1/t)}~dt.$$

- Rotating the picture with respect to $u\in S_+^1$, one gets a three-dimensional object, looking like some half pseudosphere. The corresponding volume equals, up to constants the volume of the integral we wish to illustrate.

\medskip
Let us precise that when the random walk goes frankly in some direction $u\in S^1$, then for all $v\in S_+^1$ not orthogonal to $u$, some pinching effect occurs towards the horizontal plane in the sliced picture in direction $v$, making the area (and thus the global volume) smaller.

\begin{center}
\resizebox{0.44\linewidth}{!}{\begin{picture}(0,0)%
\includegraphics{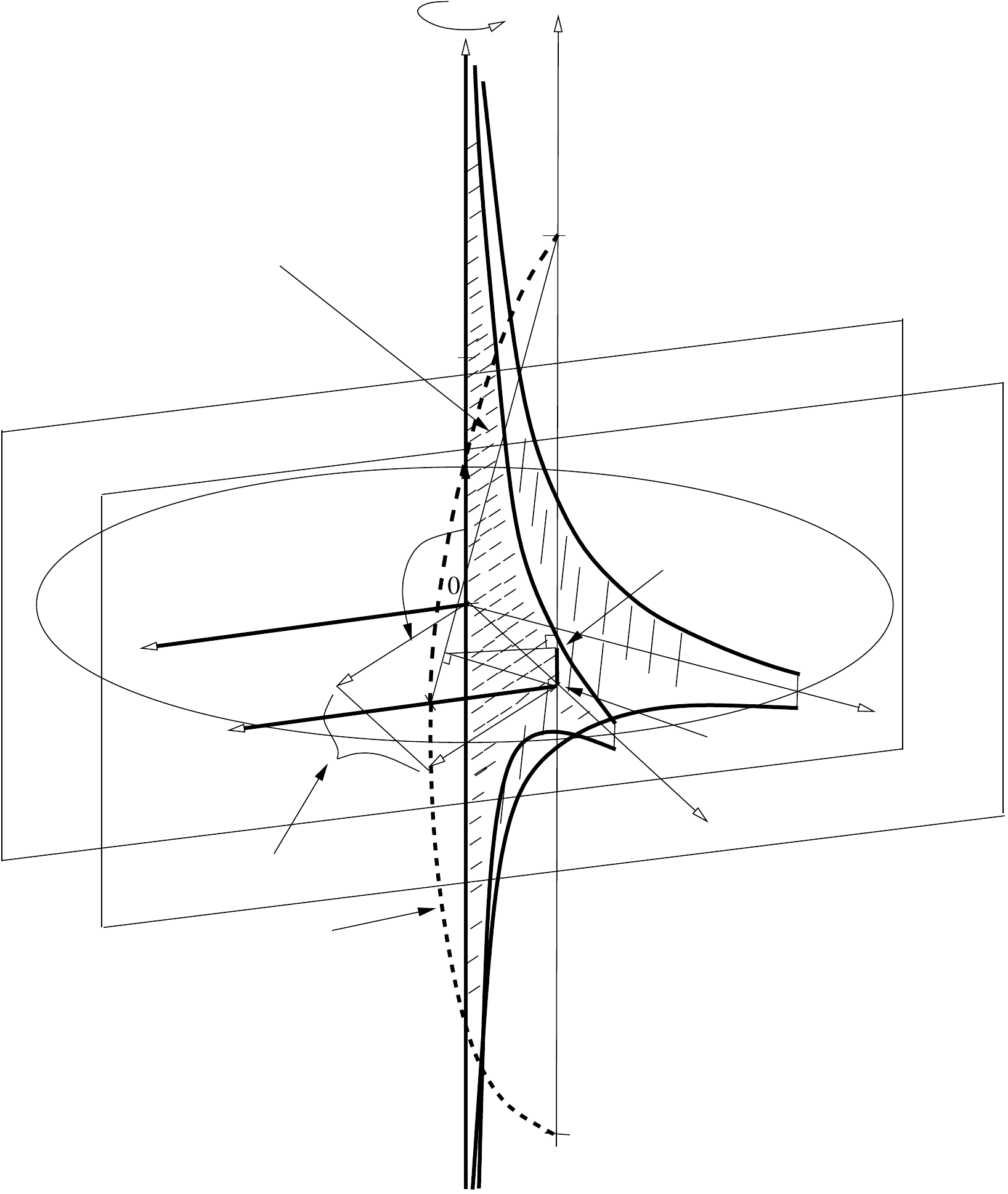}%
\end{picture}%
\setlength{\unitlength}{3947sp}%
\begingroup\makeatletter\ifx\SetFigFont\undefined%
\gdef\SetFigFont#1#2#3#4#5{%
  \reset@font\fontsize{#1}{#2pt}%
  \fontfamily{#3}\fontseries{#4}\fontshape{#5}%
  \selectfont}%
\fi\endgroup%
\begin{picture}(7839,9300)(2433,-8824)
\put(2648,-7066){\makebox(0,0)[lb]{\smash{{\SetFigFont{12}{14.4}{\rmdefault}{\mddefault}{\updefault}{\color[rgb]{0,0,0}Level line at height $t>0$ }%
}}}}
\put(3241,-1351){\makebox(0,0)[lb]{\smash{{\SetFigFont{12}{14.4}{\rmdefault}{\mddefault}{\updefault}{\color[rgb]{0,0,0}Area in direction $u$}%
}}}}
\put(5543,-3713){\makebox(0,0)[lb]{\smash{{\SetFigFont{12}{14.4}{\rmdefault}{\mddefault}{\updefault}{\color[rgb]{0,0,0}$\alpha$}%
}}}}
\put(4943,-4741){\makebox(0,0)[lb]{\smash{{\SetFigFont{12}{14.4}{\rmdefault}{\mddefault}{\updefault}{\color[rgb]{0,0,0}$\rho e^{i\alpha}$}%
}}}}
\put(2716,-6421){\makebox(0,0)[lb]{\smash{{\SetFigFont{12}{14.4}{\rmdefault}{\mddefault}{\updefault}{\color[rgb]{0,0,0}$1/\varphi_u(-\rho\sin(\alpha/2),\rho\cos(\alpha/2))$}%
}}}}
\put(2641,-7321){\makebox(0,0)[lb]{\smash{{\SetFigFont{12}{14.4}{\rmdefault}{\mddefault}{\updefault}{\color[rgb]{0,0,0}(on $0\leq\alpha\leq\pi$, $\rho>0$)}%
}}}}
\put(7988,-5348){\makebox(0,0)[lb]{\smash{{\SetFigFont{12}{14.4}{\rmdefault}{\mddefault}{\updefault}{\color[rgb]{0,0,0}$tu$}%
}}}}
\put(7658,-4019){\makebox(0,0)[lb]{\smash{{\SetFigFont{12}{14.4}{\rmdefault}{\mddefault}{\updefault}{\color[rgb]{0,0,0}$z_t$}%
}}}}
\put(5626,-2311){\makebox(0,0)[lb]{\smash{{\SetFigFont{12}{14.4}{\rmdefault}{\mddefault}{\updefault}{\color[rgb]{0,0,0}$e_3$}%
}}}}
\put(9826,-2761){\makebox(0,0)[lb]{\smash{{\SetFigFont{12}{14.4}{\rmdefault}{\mddefault}{\updefault}{\color[rgb]{0,0,0}$P_t$}%
}}}}
\put(9001,-2311){\makebox(0,0)[lb]{\smash{{\SetFigFont{12}{14.4}{\rmdefault}{\mddefault}{\updefault}{\color[rgb]{0,0,0}$P_0$}%
}}}}
\put(8678,-5123){\makebox(0,0)[lb]{\smash{{\SetFigFont{12}{14.4}{\rmdefault}{\mddefault}{\updefault}{\color[rgb]{0,0,0}$v$}%
}}}}
\put(7311,-5366){\makebox(0,0)[lb]{\smash{{\SetFigFont{12}{14.4}{\rmdefault}{\mddefault}{\updefault}{\color[rgb]{0,0,0}$u$}%
}}}}
\end{picture}%
}
\end{center}

We now discuss some consequences of the theorem.
\begin{cor}

\label{principcor}
$ $

\noindent
For the general model, a sufficient condition for transience is~:

$$\sum_{n\geq1}\int_{S_+^{d-1}}\frac{1}{(\varphi_u(n))^d}~du<+\infty.$$

\noindent
It is true under the condition $\sum_{n\geq1}\psi(n)^{-d}<+\infty$, depending only on the vertical. The latter is satisfied in the following cases~:

\smallskip

-- $d\geq3$,

-- $d=2$ and $w_+\circ v_+^{-1}(n)+w_-\circ v_-^{-1}(n)\geq (\log n)^{1+\varepsilon}$ and in particular if $p_n=q_n$, $n\in\Z$,

-- $d=1$ and $w_+\circ v_+^{-1}(n)+w_-\circ v_-^{-1}(n)\geq n(\log n)^{2+\varepsilon}$.

\end{cor}

In the antisymmetric case, an explicit criterion is available.

\begin{prop}
\label{antiscor}

$ $

\noindent
Antisymmetric case. Suppose that $m_{-n}=-m_n$ and $\rho_{-n}=\rho_n$, $n\geq0$. The random walk is transient if and only if~:
 
$$\sum_{n\geq1}\int_{S_+^{d-1}}\frac{1}{(\varphi_{u,++}(n))^d}~du<+\infty.$$

\noindent
In particular, let $m_n=-m_{-n}=c\not=0$, $n\geq1$, with $m_0=0$, and suppose that $c_1n^{\alpha}\leq \rho_n\leq c_2n^{\alpha}$, $n\geq 0$, where $\alpha\in\R$. Then~:

\smallskip

-- If $d=1$, the random walk is recurrent if and only if $\alpha\geq1$.

-- If $d=2$, the random walk is recurrent if and only if $\alpha\geq3$.

\end{prop}

There would be many other cases to consider. The Campanino-Petritis model, i.e. $d=1$, $p_n=q_n$ and $\mu_n=\delta_{\varepsilon_n}$, in the case when $\varepsilon_n=1_{n\geq0}-1_{n<0}$ corresponds to the first example with $\rho_n=1$, so $\alpha=0$, and the random walk is transient. As already indicated in \cite{cp}, the parameters are largely interior to the transience domain. Remark that in the antisymmetric case when $d=1$, taking $\mu_n=\delta_1$, $\mu_{-n}=\delta_{-1}$, for $n\geq1$, and $\rho_n\sim n^{\alpha}$, $n\geq1$, since horizontal steps are restricted to $+1$ in the North and to $-1$ in the South, the random walk (recurrent or transient) necessarily makes spirals.

\begin{center}
\resizebox{0.6\linewidth}{!}{\begin{picture}(0,0)%
\includegraphics{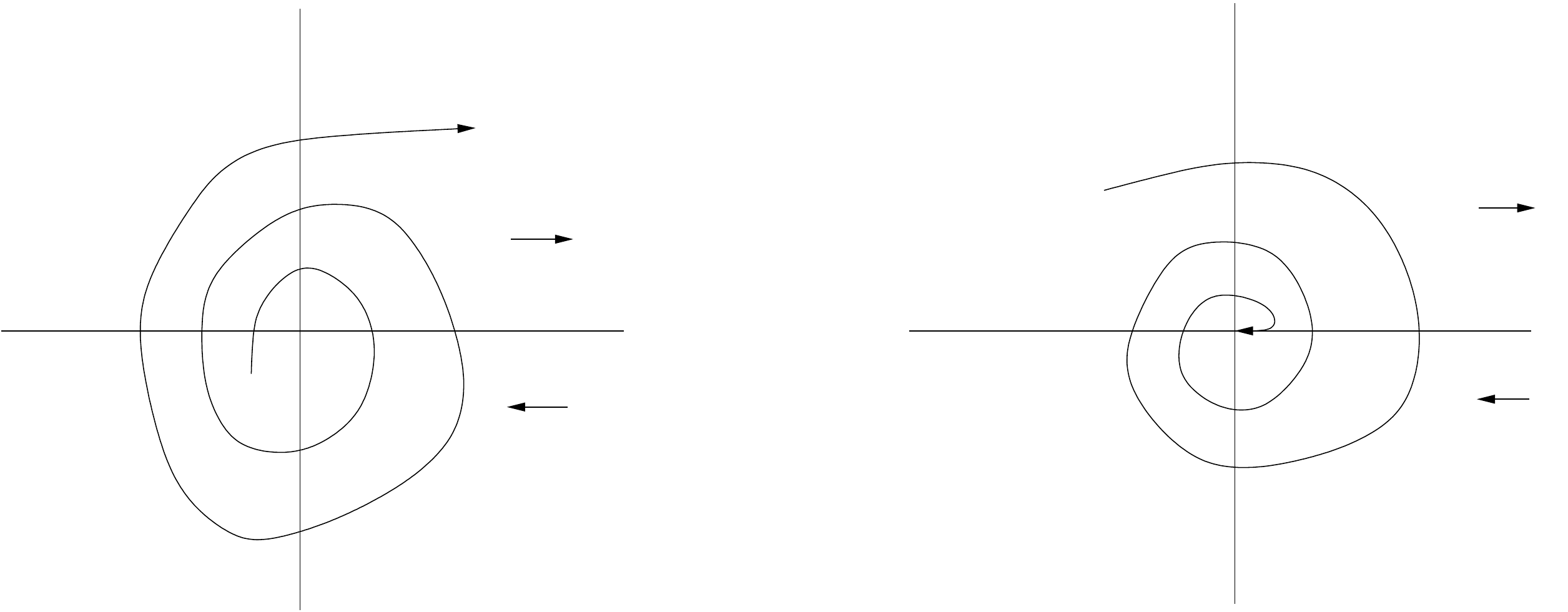}%
\end{picture}%
\setlength{\unitlength}{3947sp}%
\begingroup\makeatletter\ifx\SetFigFont\undefined%
\gdef\SetFigFont#1#2#3#4#5{%
  \reset@font\fontsize{#1}{#2pt}%
  \fontfamily{#3}\fontseries{#4}\fontshape{#5}%
  \selectfont}%
\fi\endgroup%
\begin{picture}(12042,4689)(124,-6403)
\put(7336,-2236){\makebox(0,0)[lb]{\smash{{\SetFigFont{12}{14.4}{\rmdefault}{\mddefault}{\updefault}{\color[rgb]{0,0,0}$\alpha>1$}%
}}}}
\put(4846,-3601){\makebox(0,0)[lb]{\smash{{\SetFigFont{12}{14.4}{\rmdefault}{\mddefault}{\updefault}{\color[rgb]{0,0,0}$+1$}%
}}}}
\put(4831,-4951){\makebox(0,0)[lb]{\smash{{\SetFigFont{12}{14.4}{\rmdefault}{\mddefault}{\updefault}{\color[rgb]{0,0,0}$-1$}%
}}}}
\put(12151,-3361){\makebox(0,0)[lb]{\smash{{\SetFigFont{12}{14.4}{\rmdefault}{\mddefault}{\updefault}{\color[rgb]{0,0,0}$+1$}%
}}}}
\put(12136,-4861){\makebox(0,0)[lb]{\smash{{\SetFigFont{12}{14.4}{\rmdefault}{\mddefault}{\updefault}{\color[rgb]{0,0,0}$-1$}%
}}}}
\put(346,-2251){\makebox(0,0)[lb]{\smash{{\SetFigFont{12}{14.4}{\rmdefault}{\mddefault}{\updefault}{\color[rgb]{0,0,0}$\alpha<1$}%
}}}}
\end{picture}%
}
\end{center}

In contrast with the flat case ($p_n=q_n$, $n\in\Z$), one can for this model in some sense ``suppress'' the vertical dimension for some values of the parameters. Indeed, when $\sum_{n\in\Z}(1/\rho_n)<+\infty$, the vertical component is positive recurrent, hence admits an invariant probability measure. In this sense, the random walk is then ``essentially'' $d$-dimensional. When $d=1$, this is a kind of random walk in a half-pipe. This explains the critical values of $d$ appearing in the corollary and in particular the fact that the random walk (in $\Z^{d+1}$) can be recurrent when $d=2$. 

\begin{prop}
\label{halfpipe}

$ $

\noindent
Suppose that $\sum_{n\in\Z}(1/\rho_n)<+\infty$ and $d=1$. 

\smallskip
\noindent
i) If $\sum_{n\in\Z}(m_n/\rho_n)\not=0$, then the random walk is transient.

\medskip
\noindent
ii) If $\sum_{n\in\Z}(m_n/\rho_n)=0$ and $p_{-n}=q_n$, $r_{-n}=r_n$, $\mu_n=\mu_{-n}$, for $n\geq0$, then the random walk is recurrent.
\end{prop}

Hence, when $\sum_{n\in\Z}(1/\rho_n)<+\infty$, the finiteness condition for transience is replaced by some non-zero condition. Fixing $d=1$ and, breaking momentarily the assumptions, suppose that brutally $\rho_1=\rho_{-1}=+\infty$ (giving $\rho_k=+\infty$, $k\not=0$). One then recovers that the condition $m_0\not=0$ is necessary and sufficient for transience, which is a standard result for one-dimensional i.i.d random walk with integrable step. 

\medskip
\noindent
Notice in such a model the important role of a single hyperplane, as the latter can modify the asymptotics. This is not true if $\sum_{n\in\Z}1/\rho_n$, for example if $p_n=q_n$, $n\in\Z$, and $d=1$, where changing one line did not modify the asymptotics (see the introduction in \cite{jb1}).

 \medskip
 We give an application to a random walk in a half-pipe, with independent level lines.

\begin{cor} 

$ $

\noindent
Let $d=1$ and $\sum_{n\in\Z}(1/\rho_n)<+\infty$. Suppose that the $(m_n)_{n\in\Z}$ are a typical realization of some independent uniformly bounded random variables, at least one having a density. Then, almost-surely, the associated random walk is transient. 
\end{cor}

\medskip
Indeed, it is clear from the hypotheses that the random variable $\omega\longmapsto\sum_{n\in\Z}(m_n(\omega)/\rho_n)$ admits a density, so equals $0$ with zero probability. We next apply the result of the previous proposition.

\begin{center}
\resizebox{0.40\linewidth}{!}{\begin{picture}(0,0)%
\includegraphics{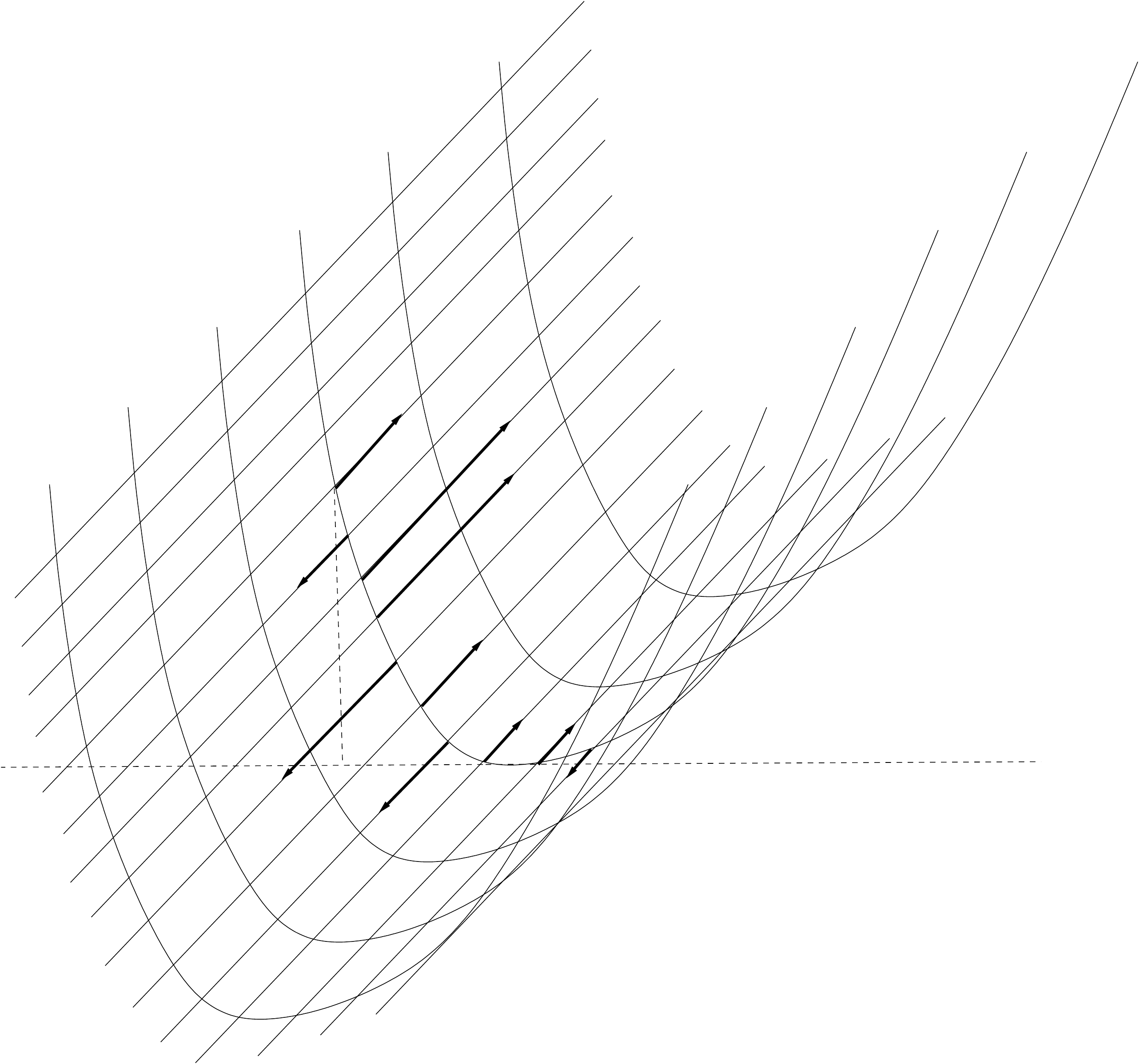}%
\end{picture}%
\setlength{\unitlength}{3947sp}%
\begingroup\makeatletter\ifx\SetFigFont\undefined%
\gdef\SetFigFont#1#2#3#4#5{%
  \reset@font\fontsize{#1}{#2pt}%
  \fontfamily{#3}\fontseries{#4}\fontshape{#5}%
  \selectfont}%
\fi\endgroup%
\begin{picture}(12305,11499)(664,-9223)
\put(3721,-3016){\makebox(0,0)[lb]{\smash{{\SetFigFont{12}{14.4}{\rmdefault}{\mddefault}{\updefault}{\color[rgb]{0,0,0}$(0,n)$}%
}}}}
\put(5401,-6136){\makebox(0,0)[lb]{\smash{{\SetFigFont{12}{14.4}{\rmdefault}{\mddefault}{\updefault}{\color[rgb]{0,0,0}$0$}%
}}}}
\put(4021,-4741){\makebox(0,0)[lb]{\smash{{\SetFigFont{12}{14.4}{\rmdefault}{\mddefault}{\updefault}{\color[rgb]{0,0,0}$\rho_n$}%
}}}}
\put(4561,-2161){\makebox(0,0)[lb]{\smash{{\SetFigFont{12}{14.4}{\rmdefault}{\mddefault}{\updefault}{\color[rgb]{0,0,0}$m_n$}%
}}}}
\end{picture}%
}
\end{center}

\noindent
In this picture in $\R^3$ of the $\Z^2$-half-pipe, we have drawn the points $(k,l)\in\Z^2$ at height $\rho_l$. The quantity $\rho_l$ can be considered as the ``level of the sea" at $(k,l)$. The borders of the half-pipe are very steep due to the condition $\sum_{n\in\Z}(1/\rho_n)<+\infty$.

\section{Preliminaries}
\subsection{Sleszynski-Pringsheim continued fractions}
Formally, a general finite continued fraction is written as follows~:

$$[(c_1,d_1); (c_2,d_2);\cdots;(c_n,d_n)]=\cfrac{c_1}{d_1+\cfrac{c_2}{d_2+\cfrac{\cdots}{\cdots+ \cfrac{c_n}{d_n}}}}.$$

\smallskip
\noindent
We shall consider finite continued fractions corresponding to the application to some $z_0\in\C$ in the unit disk  of functions of the form $z\longmapsto c/(d+z)$, with complex numbers $c\not=0$ and $d$ so that $|c|+1\leq |d|$, hence preserving the closed unit disk. Such finite continued fractions are usually called finite Sleszynski-Pringsheim (SP) continued fractions. 

\medskip
Infinite SP-continued fractions, written $[(c_1,d_1); (c_2,d_2);\cdots]$, also converge, by the Sleszynski-Pringsheim theorem (see \cite{lw}). We will reproduce the arguments proving this result.

\medskip
For $n\geq0$, the finite continued fraction $[(c_1,d_1); (c_2,d_2);\cdots;(c_n,d_n)]$ can be reduced as a fraction $A_n/B_n$, where the $(A_n)$ and $(B_n)$ satisfy the recursive relations~:

$$
\left\{{\begin{array}{c}
A_n=d_nA_{n-1}+c_nA_{n-2},~n\geq 1,~A_{-1}=1,~A_0=0,\\
\\
B_n=d_nB_{n-1}+c_nB_{n-2},~n\geq 1,~B_{-1}=0,~B_0=1.\\
\end{array}}\right.$$

\smallskip
\noindent
In our setting it will be directly checked that $B_n$ is never zero for $n\geq0$. We require the following classical determinant. For $n\geq 1$~:

\begin{eqnarray}
A_nB_{n-1}-A_{n-1}B_n&=&(-c_n)(A_{n-1}B_{n-2}-A_{n-2}B_{n-1})=\cdots\nonumber\\
&=&(-1)^nc_1\cdots c_n(A_{0}B_{-1}-A_{-1}B_{0})=(-1)^{n+1}c_1\cdots c_n.\nonumber
\end{eqnarray}

\medskip
\noindent
This gives the following representation as a series~:

\begin{equation}
\label{series}\displaystyle[(c_1,d_1); (c_2,d_2);\cdots;(c_n,d_n)]=\frac{A_n}{B_n}=\sum_{k=1}^n\({\frac{A_k}{B_k}-\frac{A_{k-1}}{B_{k-1}}}\)=\sum_{k=1}^n\frac{(-1)^{k+1}c_1\cdots c_k}{B_kB_{k-1}}.
\end{equation}

\medskip
We now focus on a particular class of SP-continued fractions that will appear frequently.

 \begin{lemme}
 \label{cf}
 
 $ $
 
 \noindent
 Assume that $\lim_{n\rightarrow+\infty}v_+(n)=+\infty$.
 $ $
 
\begin{enumerate}
\item Let $(\gamma_n)_{n\geq 1}$ and $(\gamma'_n)_{n\geq 1}$ be sequences of complex numbers with $0<|\gamma_n|\leq1$, $|\gamma'_n|\leq 1$. Then~:

$$[(a_1,b_1/\gamma_1); (-a_2,b_2/\gamma_2);\cdots;(-a_{n-1},b_{n-1}/\gamma_{n-1});(-a_n,b_n/\gamma_n-\gamma'_n)]$$

\noindent
is well-defined. It converges to $[(a_1,b_1/\gamma_1); (-a_2,b_2/\gamma_2);\cdots;(-a_n,b_n/\gamma_n);\cdots]$, as $n\rightarrow+\infty$, an infinite SP-continued fraction. The latter is the limit of $A_n/B_n$, as $n\rightarrow+\infty$, where~:
 
\begin{equation}
\label{syst}
\left\{{\begin{array}{c}
\displaystyle A_n=\frac{b_n}{\gamma_n}A_{n-1}-a_nA_{n-2},~n\geq 2,~A_{-1}=1,~A_0=0,~A_1=a_1,\\
\\
\displaystyle B_n=\frac{b_n}{\gamma_n}B_{n-1}-a_nB_{n-2},~n\geq 2,~B_{-1}=0,~B_0=1,~B_1=b_1/\gamma_1.\\
\end{array}}\right.
\end{equation}
\item Set $v_+(-1)=0$. The solutions $(B_n)$ of \eqref{syst} check~:

$$|B_n|-|B_{n-1}|\geq a_n(|B_{n-1}|-|B_{n-2}|),~n\geq 1.$$

\noindent
As a result, $|B_n|\geq v_+(n),~n\geq-1$. If the $0<\gamma_n\leq1$ are real, then $B_n>B_{n-1}>\cdots> B_{-1}=0$. When $\gamma_n=1$, $n\geq1$, then $B_n=v_+(n)$, $n\geq -1$, as well as $A_n=v_+(n)-1$, $n\geq 0$.
\item In \eqref{syst}, $n\longmapsto |B_n|/v_+(n)$, $n\geq0$, is non-decreasing. Also, for $n\geq1$~:

$$\sum_{k>n}\frac{\rho_k}{|B_kB_{k-1}|}\leq \frac{v_+(n)}{|B_n|^2}\leq \frac{1}{|B_n|}.$$

\end{enumerate}
\end{lemme} 
 
\noindent
\textit{Proof of the lemma~:}

\noindent
The solutions of \eqref{syst} check $|B_n|\geq b_n|B_{n-1}|-a_n|B_{n-2}|$, $n\geq1$. Hence~:

$$|B_n|-|B_{n-1}|\geq a_n(|B_{n-1}|-|B_{n-2}|),~n\geq1.$$

\smallskip
\noindent
When iterating, $|B_n|-|B_{n-1}|\geq \rho_n$. Thus $|B_n|\geq v_+(n).$ 

\medskip
In point $1.$, the finite continued fraction is well-defined because $a_k\not=0$, $|b_k/\gamma_k|-a_k\geq b_k-a_k=1$ and $|\gamma'_n|\leq 1$. We obtain from \eqref{series}~:

\begin{equation}
\label{cvv}
[(a_1,b_1/\gamma_1); (-a_2,b_2/\gamma_2);\cdots;(-a_n,b_n/\gamma_n-\gamma'_n)]=\sum_{k=1}^{n-1}\frac{\rho_k}{B_kB_{k-1}}+\frac{\rho_n}{B_{n-1}\tilde{B}_n},
\end{equation}

\noindent
where $\tilde{B}_n=(b_n/\gamma_n-\gamma'_n)B_{n-1}-a_nB_{n-2}$. We get~:

$$|\tilde{B}_n|\geq (b_n-1)|B_{n-1}|-a_n|B_{n-2}|\geq a_n(|B_{n-1}|-|B_{n-2}|)\geq a_n\rho_{n-1}=\rho_n.$$

\medskip
\noindent
In \eqref{cvv} the first term in the right-hand side is absolutely convergent, because, as $|B_k|\geq v_+(k)$~:

$$\sum_{k\geq1}\frac{\rho_k}{|B_kB_{k-1}|}\leq\sum_{k\geq1}\frac{v_+(k)-v_+(k-1)}{v_+(k)v_+(k-1)}=\sum_{k\geq1}\({\frac{1}{v_+(k-1)}-\frac{1}{v_+(k)}}\)=1.$$
 
 \noindent
As $|B_{n-1}|\rightarrow+\infty$, we conclude that the right-hand side in \eqref{cvv} converges to $\sum_{k\geq 1}\rho_k/(B_kB_{k-1})$.

\medskip
When the $\gamma_n$ are real, write $B_n-B_{n-1}=\frac{(1-\gamma_n)b_n}{\gamma_n}B_{n-1}+a_n(B_{n-1}-B_{n-2})$. The condition ``$B_n> B_{n-1}\geq 0$'' is then transmitted recursively. If $\gamma_n=1$, then $B_n-B_{n-1}=a_n(B_{n-1}-B_{n-2})$, giving $B_n=v_+(n)$, $n\geq0$. Similarly, $A_n=v_+(n)-1$, $n\geq0$, as~:

$$A_n-A_{n-1}=a_n(A_{n-1}-A_{n-2})=\cdots=a_n\cdots a_2(A_1-A_0)=\rho_n\mbox{ and }A_0=0.$$ 

\medskip
For the last point, we first show that $n\longmapsto |B_n|/v_+(n)$, $n\geq0$, is non-decreasing. We will require it in the equivalent form $|B_n|/(|B_{n+1}|-|B_n|)\leq v_+(n)/(v_+(n+1)-v_+(n))$. Write~:

\begin{eqnarray}
v_+(n)|B_{n+1}|-v_+(n+1)|B_n|&\geq& v_+(n)\({b_{n+1}|B_n|-a_{n+1}|B_{n-1}|}\)-v_+(n+1)|B_n|\nonumber\\
&\geq&a_{n+1}\({v_+(n-1)|B_n|-v_+(n)|B_{n-1}|}\)\nonumber\\
&\geq&\cdots\geq \rho_{n+1}\({|B_0|v_+(-1)-v_+(0)|B_{-1}|}\)=0.\nonumber\end{eqnarray}

\noindent
Finally, using the previous results, for $n\geq1$~:

\begin{eqnarray}
\sum_{k>n}\frac{\rho_k}{|B_kB_{k-1}|}&=&\sum_{k>n}\rho_k\({\frac{1}{|B_{k-1}|}-\frac{1}{|B_k|}}\)\frac{1}{|B_k|-|B_{k-1}|}\nonumber\\
&\leq&\sum_{k>n}\rho_k\({\frac{1}{|B_{k-1}|}-\frac{1}{|B_k|}}\)\frac{1}{a_k\cdots a_{n+2}(|B_{n+1}|-|B_{n}|)}\nonumber\\
&\leq&\frac{\rho_{n+1}}{|B_{n+1}|-|B_n|}\sum_{k>n}\({\frac{1}{|B_{k-1}|}-\frac{1}{|B_k|}}\)\leq\frac{v_+(n+1)-v_+(n)}{|B_{n+1}|-|B_n|}\frac{1}{|B_{n}|}\leq \frac{v_+(n)}{|B_n|^2}.\nonumber\end{eqnarray}

\noindent
This completes the proof of the lemma.\fin

\subsection{Asymptotical behavior of the vertical component}
The question of the recurrence/transience of the vertical component of the random walk is classical. Indeed the vertical component restricted to the subsequence of vertical movements is the random walk on $\Z$ with transition probabilities $\P_{n,n-1}=q_n'$ and $\P_{n,n+1}=p'_n$, $n\in\Z$.

\begin{lemme}

$ $

\noindent
The Markov chain on $\Z$ so that $\P_{n,n+1}=p'_n$ and $\P_{n,n-1}=q'_n$, $n\in\Z$, is recurrent if and only if $\lim_{n\rightarrow+\infty}v_+(n)=+\infty$ and $\lim_{n\rightarrow+\infty}v_-(n)=+\infty$.
\end{lemme}

\noindent
{\it Proof of the lemma~:}

\noindent
Fix $N>1$ and let $f(k)=\P_k(\mbox{exit }[0,N]\mbox{ on the left side})$,~$0\leq k\leq N$. The Markov property implies that $k\longmapsto f(k)$ is harmonic in the interior of this interval. Precisely, for $1\leq k\leq N-1$~:

$$f(k)=p'_kf(k+1)+q'_kf(k-1).$$

\smallskip
\noindent
Let $g(k)=f(k)-f(k-1)$. We obtain $g(k)=(p_k/q_k)g(k+1)$ and therefore $g(k)=\rho_{k-1}g(1)$, $1\leq k\leq N$. As a result~:

$$-1=\sum_{k=1}^Ng(k)=-\P_1(\mbox{exit }[0,N]\mbox{ at }N)\sum_{1\leq k\leq N}\rho_{k-1}.$$

\noindent
Hence $\P_1(\mbox{reach }0)=1\Leftrightarrow\lim_{n\rightarrow+\infty}v_+(n)=+\infty$. Idem $\P_{-1}(\mbox{reach }0)=1\Leftrightarrow\lim_{n\rightarrow+\infty}v_-(n)=+\infty$. This furnishes the desired result.\fin

\bigskip
The previous criterion can be reformulated using trees. Let us say that a random variable $X$ has the geometrical law ${\cal G}(p)$, $0<p<1$, if $\P(X=n)=p^n(1-p)$, $n\geq 0$. 

\begin{lemme}

$ $

\noindent
Consider the Galton-Watson tree $(Z^+_n)_{n\geq1}$ with $Z^+_1=1$ and, independently, the law of the number of children at level $n+1$ of an individual at level $n\geq1$ is ${\cal G}(p'_n)$. Then this tree is finite almost-surely if and only if $\lim_{n\rightarrow+\infty}v_+(n)=+\infty$.
\end{lemme}

\noindent
{\it Proof of the lemma~:}

\noindent
As usual, since $\{Z^+_n=0\}\subset\{Z^+_{n+1}=0\}$, the almost-sure finiteness is equivalent to  $\P(Z^+_n=0)\rightarrow1$. Fix $0<s<1$ and recall that $\E(s^{Z^+_n})-s\leq\P(Z^+_n=0)\leq \E(s^{Z^+_n})$. Taking $n\geq 2$~:

$$\E\left({s^{Z^+_n}}\right)=\E\left[{\left({\frac{1-p'_{n-1}}{1-sp'_{n-1}}}\right)^{Z^+_{n-1}}}\right]=\E\left[{\left({\frac{a_{n-1}}{b_{n-1}-s}}\right)^{Z^+_{n-1}}}\right].$$

\noindent
Iterating (using $a_{n-1}/(b_{n-1}-s)$ in place of $s$), we obtain the following SP-continued fraction~:

$$\E\left({s^{Z^+_n}}\right)=\[{(a_1,b_1);(-a_2,b_2);\cdots;(-a_{n-2},b_{n-2});(-a_{n-1},b_{n-1}-s)}\].$$

\smallskip
\noindent
This corresponds to $\gamma_k=1$ and $\gamma'_n=s$ in lemma \ref{cf}. From lemma \ref{cf} and relation \eqref{cvv}~:

$$\E\left({s^{Z^+_n}}\right)=\frac{v_+(n-2)-1}{v_+(n-2)}+\frac{\rho_{n-1}}{v_+(n-1)\tilde{B}_{n-1}},$$

\noindent
with $\tilde{B}_{n-1}=(b_{n-1}-s)v_+(n-2)-a_{n-1}v_+(n-3)$, so that $\tilde{B}_{n-1}\geq\rho_{n-1}$ and $\tilde{B}_{n-1}\geq (1-s)v_+(n-2)$. If $v_+(n)\rightarrow+\infty$, then $\E(s^{Z^+_n})\rightarrow1$ uniformly in $0<s<1$, giving $\P(Z^+_n=0)\rightarrow1$. If $v_+(n)\rightarrow_{n\rightarrow+\infty}b\in(0,+\infty)$, then $\rho_n\rightarrow0$ and for fixed $0<s<1$ we have $\liminf_n\tilde{B}_{n-1}\geq(1-s)b>0$, so that $\E(s^{Z^+_n})$ tends to $(b-1)/b<1$, giving $\lim_n \P(Z^+_n=0)=(b-1)/b$.\fin

\bigskip
\noindent
\begin{remark}
There is naturally a symmetric result for the Southern direction of the vertical component. One introduces, with decreasing indices $n\leq-1$, the Galton-Watson tree $(Z^-_n)_{n\leq-1}$ with $Z^-_{-1}=1$ such that, independently, the law of the number of children at level $n-1$ of an individual at level $n$ is ${\cal G}(q'_{n})$. The tree is almost-surely finite if and only if $\lim_{k\rightarrow+\infty}v_-(k)=+\infty$. \end{remark}

\section{Reduction to an i.i.d. random walk in $\Z^d$}
For the rest of the article we therefore suppose the vertical component recurrent. Equivalently, from the previous section, this means $\lim_{n\rightarrow+\infty}v_{\pm}(n)=+\infty$. Just observe that if for example $v_+$ is bounded by some $v_+(\infty)<\infty$, then $+\infty=\varphi_u(n)=\varphi_{u,+}(n)=\psi(n)$ for $n>v_+(\infty)$, so the reversed functions are bounded quantities and the integral involved in the theorem is finite. The same happens if $v_-$ is bounded.

\medskip
We can now introduce the random times $0=\sigma_0<\tau_0<\sigma_1<\tau_1<\cdots$, where~:

$$\tau_k=\min\{n>\sigma_k~|~S^{2}_n\not=0\},~\sigma_{k+1}=\{n>\tau_k~|~S^{2}_n=0\}.$$

\medskip
\noindent
Introduce the $\Z^d$-displacement $D_n=S^1_{\sigma_n}-S^1_{\sigma_{n-1}}$. As the environment is invariant under $\Z^d$-translations, the $(D_n)_{n\geq 1}$ are globally independent and identically distributed. The following lemma is essentially contained in \cite{cp}.

\begin{lemme}
\label{petri}

$ $

\noindent
Let $T_0=0$ and $T_n=D_1+\cdots+D_n$, $n\geq1$. The random walk $(S_n)_{n\geq 0}$ is recurrent is $\Z^{d+1}$ if and only if $(T_n)_{n\geq 0}$ is recurrent in $\Z^d$.
\end{lemme}

\noindent
{\it Proof of the lemma~:}

\noindent
If $(T_n)_{n\geq 0}$ is recurrent in $\Z^d$, then $(S_n)$ is recurrent in $\Z^{d+1}$, as $S_{\sigma_n}=(T_n,0)$. In case of transience of $(T_n)$, using again the invariance of the environment under $\Z^d$-translations, we have~:

$$\exists C,~\forall x\in\Z^d,~\sum_{n\geq 1}\P\({T_n=x}\)\leq C.$$

\noindent
Let $\Gamma\sim{\cal G}(r_0)$ and $\xi_k\sim\mu_0$, for $k\geq 1$, so that $((\xi_k)_{k\geq 1},\Gamma)$ are globally independent and also from the sequence $(T_n)$. Remark that $(S_l^1)_{l\in [\sigma_k,\tau_k)}$ and $(T_k+\sum_{1\leq m\leq l}\xi_m)_{0\leq l\leq \Gamma}$ have the same law. Introduce the real random variable~:

$$H=\sum_{1\leq k\leq \Gamma}\|\xi_k\|.$$

\noindent
Observe now that $S_n$ can be $0$ only for $n$ in some $[\sigma_k,\tau_k)$ and that~:

$$\P(\exists n\in[\sigma_k,\tau_k),~S_n=0)\leq\P(H\geq \|T_k\|).$$

\noindent
This provides~:

\begin{eqnarray}
\sum_{k\geq 1}\P(\exists n\in[\sigma_k,\tau_k),~S_n=0)\leq\sum_{k\geq 1}\P(H\geq\|T_k\|)&\leq&\sum_{x\in\Z^d}\sum_{k\geq 1}\P(T_k=x)\P\({H\geq \|x\|}\)\nonumber\\
&\leq&C\sum_{x\in\Z^d}\P(H\geq \|x\|)\leq C'\E(H^d).\nonumber\end{eqnarray}

\noindent
Finally, this gives~:

\begin{eqnarray}
\E(H^d)=\sum_{n\geq 0}\P(\Gamma=n)\E\[{\({\sum_{1\leq k\leq n}\|\xi_k\|}\)^d}\]&\leq& (1-r_0)\sum_{n\geq 0}r_0^nn^{d-1}\E\({\sum_{1\leq k\leq n}\|\xi_k\|^d}\)\nonumber\\
&\leq&(1-r_0)\sum_{n\geq 0}r_0^nn^{d}\E(\|\xi_1\|^d)<\infty.\nonumber
\end{eqnarray}

\noindent
By the Borel-Cantelli lemma, $(S_n)$ is transient. This completes the proof of the lemma. \fin

\bigskip
This reduces the problem of the recurrence of $(S_n)$ to that of $(T_n)$. Set~:

$$D=D_1\mbox{ and }\chi_D(t)=\E(e^{it.D}),~t\in\R^d.$$

\medskip
\noindent
We shall use the following theorem, the strong form of the Chung-Fuchs recurrence criterion, giving an analytical recurrence criterion for a $i.i.d.$ random walk in $\Z^d$. See Spitzer \cite{spitz}. Recall that $S_+^{d-1}$ denotes the half unit sphere and let $B_d(0,\eta)$ be the ball of center $0$ and radius $\eta>0$ in $\R^d$.

\begin{thm}
\label{spitz}

$ $ 

\noindent
Suppose that the subgroup of $(\Z^d,+)$ generated by the support of the law of $D$ is $\Z^d$. Then the random walk $(T_n)_{n\geq 0}$ is transient if and only if for some $\eta>0$~:
\begin{equation}
\label{recr}
\int_{B_d(0,\eta)}\mbox{Re}\({\frac{1}{1-\chi_D(x)}}\)~dx<+\infty.
\end{equation}
\end{thm}

\noindent
Notice that one can restrict the integral to the half unit ball $S_+^{d-1}.]0,\eta[$. Forgetting the multiplicative constant coming from the change of variables in polar coordinates, we next decompose the integral in the form~:

$$\int_{S_+^{d-1}\times]0,\eta[}\mbox{Re}\({\frac{1}{1-\chi_D(ut)}}\)~t^{d-1}dudt\mbox{, with }(u,t)\in S_+^{d-1}\times]0,\eta[.$$

\smallskip
\noindent
From our assumptions, the subgroup $G_D$ of $(\Z^d,+)$ generated by the support of the law of $D$ is $d$-dimensional. Observe that $(T_n)$ lives in $G_D$ and recall that $G_D$ admits a basis over $\Z$. Reparametrizing $G_D$ corresponds to making a linear change of variables in \eqref{recr}. The properties of dominated variations shown below in lemma \ref{young} imply that we can assume that $G_D=\Z^d$ from the beginning. This is what we do in the sequel. 

\medskip
The only singularity of $1/(1-\chi_D)$ in $\R^d/\Z^d$ is now $0$. Fixing $0<\eta<1/2$ small enough, we take $u\in S_+^{d-1}$ and $0<t<\eta$. 

\subsection{Local time and contour of a Galton-Watson tree}
For $u\in S_+^{d-1}$ we study the behavior near $0^+$ of $t\longmapsto\chi_D(ut)$. Let us introduce the one-dimensional random walk $(Y_n)_{n\geq 0}$ on $\Z$ such that $Y_0=0$ and $\P_{n,n-1}=q'_n$ and $\P_{n,n+1}=p'_n$, for $n\in\Z$. This is $(S^2_n)_{n\geq 0}$ restricted to the sequence of vertical jumps.

\medskip
\noindent
Let $\sigma=\min\{k\geq 1~|~Y_k=0\}$ be the return time to $0$. Grouping in packets the successive $\Z^d$-steps of the random walk, observe that $D$ can be written as~:

$$D=\sum_{k=0}^{\sigma-1}\({\sum_{m=1}^{\Gamma_k}\xi_m^{(k)}}\),$$

\noindent
where, conditionally on the $(Y_l)_{l\geq 0}$, the $((\xi_m^{(k)})_{m\geq 1,k\geq 0},(\Gamma_k)_{k\geq 0})$ are independent with $\xi_m^{(k)}\sim\mu_{Y_k}$ and $\Gamma_k\sim({\cal G}(r_{Y_k}))$, for all $k\geq 0$. To detail $\chi_D$, define for $n\in\Z$~:

\begin{equation}
\label{varphin}
\varphi_n(ut)=\E\({\exp\({itu.\sum_{m=1}^{\Gamma}\xi_m}\)}\),~t\in\R,
\end{equation}

\smallskip
\noindent
with random variables $\Gamma\sim {\cal G}(r_n)$ and $\xi_m\sim \mu_n$, for $m\geq1$, all being independent. Conditioning on the $(Y_l)_{l\geq 0}$, we obtain the equality~:

$$\chi_D(ut)=\E\({\prod_{k=0}^{\sigma-1}\varphi_{Y_k}(ut)}\)=\varphi_0(ut)\E\({\prod_{k=1}^{\sigma-1}\varphi_{Y_k}(ut)}\).$$

\noindent
\smallskip
The only remaining alea is that of the $(Y_l)_{l\geq 0}$. Introduce the conditional expectations~:

$$\E^+(.)=\E(.~|~Y_1=1)\mbox{ and }\E^-(.)=\E(.~|~Y_1=-1).$$

\medskip
\noindent
Setting $\chi_D^{\pm}(ut)=\E^{\pm}\({\prod_{k=1}^{\sigma-1}\varphi_{Y_k}(ut)}\)$, this leads to~:

\begin{eqnarray}
\label{gw}
\chi_D(ut)&=&\varphi_0(ut)(p'_0\chi_D^+(ut)+q'_0\chi^-_D(ut)).\end{eqnarray}

\smallskip
We next restrict the analysis to $\chi_D^+$, the case of $\chi_D^-$ being symmetric. Introducing the local times $N_n=\#\{1\leq k\leq\sigma-1,~Y_k=n\}$, $n\geq 1$, we obtain~:

$$\chi_D^+(ut)=\E^+\({\prod_{n\geq1}(\varphi_{n}(ut))^{N_n}}\).$$

\smallskip
\noindent
The alea now is on the $(N_n)_{n\geq 1}$. To describe these local times, one classically introduces (cf \cite{legall} for instance) the Galton-Watson tree $(Z^+_n)_{n\geq1}$ with $Z^+_1=1$ such that, independently, the law of the number of children at level $n+1$ of an individual at level $n$ is ${\cal G}(p'_n)$. This tree is almost-surely finite, from the hypothesis $\lim_{n\rightarrow+\infty}v_+(n)=+\infty$. 

\begin{center}
\resizebox{0.7\linewidth}{!}{\begin{picture}(0,0)%
\includegraphics{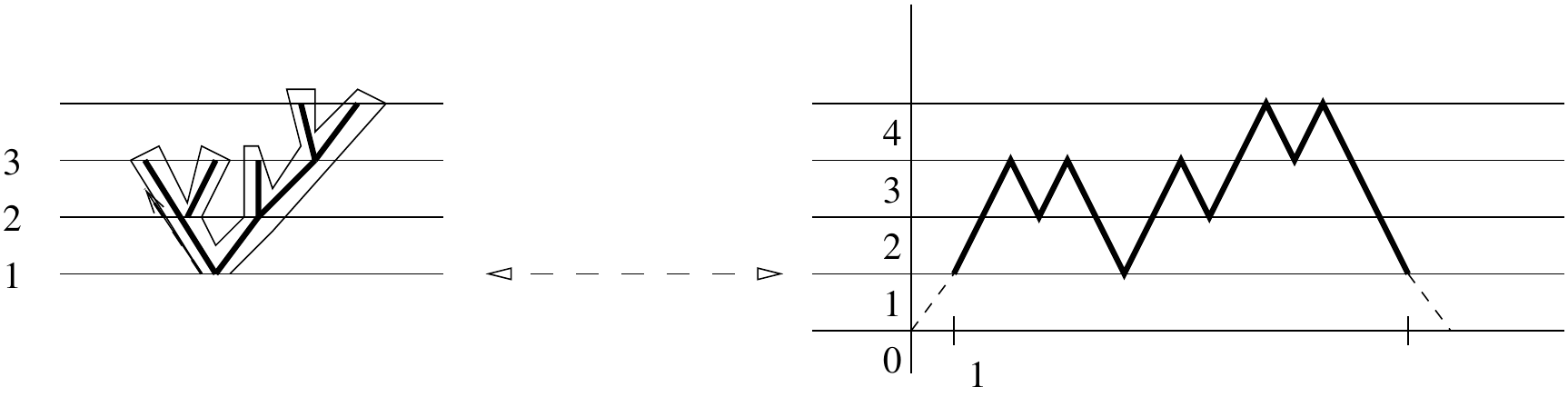}%
\end{picture}%
\setlength{\unitlength}{3947sp}%
\begingroup\makeatletter\ifx\SetFigFont\undefined%
\gdef\SetFigFont#1#2#3#4#5{%
  \reset@font\fontsize{#1}{#2pt}%
  \fontfamily{#3}\fontseries{#4}\fontshape{#5}%
  \selectfont}%
\fi\endgroup%
\begin{picture}(8277,2190)(1936,-5014)
\put(1951,-3436){\makebox(0,0)[lb]{\smash{{\SetFigFont{12}{14.4}{\rmdefault}{\mddefault}{\updefault}{\color[rgb]{0,0,0}4}%
}}}}
\put(9076,-4936){\makebox(0,0)[lb]{\smash{{\SetFigFont{12}{14.4}{\rmdefault}{\mddefault}{\updefault}{\color[rgb]{0,0,0}$\sigma-1$}%
}}}}
\end{picture}%
}
\end{center}

\noindent
As shown on the left-hand side of the picture, we make the contour process of the tree, starting from the root of the tree and turning clockwise. We associate to each ascending/descending movement a $+1/-1$ step. This gives the picture on the right-hand side, where we recover a positive excursion of the random walk $(Y_n)$ in the time interval $[1,\sigma-1]$. 

\medskip
\noindent
Observe that the total number of visits of the random walk at level $n\geq 1$ is $N_n=Z^+_n+Z^+_{n+1}$. This furnishes~:

$$\prod_{n\geq1}(\varphi_{n}(ut))^{N_n}=\prod_{n\geq1}(\varphi_{n}(ut))^{Z^+_n+Z^+_{n+1}}=\varphi_1(ut)\prod_{n\geq 1}[\varphi_n(ut)\varphi_{n+1}(ut)]^{Z^+_{n+1}}.$$

\noindent
Finally~:

$$\chi_D^+(ut)=\varphi_1(ut)\E^+\({\prod_{n\geq 1}[\varphi_n(ut)\varphi_{n+1}(ut)]^{Z^+_{n+1}}}\).$$

\subsection{Development of $\chi_D^+$ in SP-continued fraction}

We now express $\chi_D^+$ as a SP-continued fraction. For $N\geq 1$ set~:

\begin{equation}
\label{chiinit}
\chi_D^{+,N}(ut)=\varphi_1(ut)\E^+\({\prod_{n=1}^N[\varphi_n(ut)\varphi_{n+1}(ut)]^{Z^+_{n+1}}}\).
\end{equation}

\smallskip
\noindent
Let $(R_k^{(n)})_{n\geq 1,k\geq 1}$ be independent random variables such that $R_k^{(n)}\sim{\cal G}(p'_n)$. Then $(Z^+_n)_{n\geq 1}$ admits the following classical description~:

$$Z^+_1=1,~Z^+_{n+1}=\sum_{k=1}^{Z^+_n}R_k^{(n)},~n\geq 1.$$

\smallskip
\noindent
Recall that the generating function of ${\cal G}(p'_n)$ is $s\longmapsto q'_n/(1-p'_ns)=a_n/(b_n-s)$, $0\leq s\leq 1$. Using conditioning on the first step, this allows to write~:

\begin{eqnarray}\chi_D^{+,N}(ut)&=&\varphi_1(ut)\E^+\({\prod_{n=1}^{N-1}[\varphi_n(ut)\varphi_{n+1}(ut)]^{Z^+_{n+1}}(\varphi_N(ut)\varphi_{N+1}(ut))^{Z^+_{N+1}}}\)\nonumber\\
&=&\varphi_1(ut)\E^+\({\prod_{n=1}^{N-1}[\varphi_n(ut)\varphi_{n+1}(ut)]^{Z^+_{n+1}}\({\frac{a_N}{b_N-\varphi_N(ut)\varphi_{N+1}(ut)}}\)^{Z^+_N}}\)\nonumber\\
&=&\varphi_1(ut)\E^+\({\prod_{n=1}^{N-2}[\varphi_n(ut)\varphi_{n+1}(ut)]^{Z^+_{n+1}}\({\frac{a_N\varphi_{N-1}(ut)\varphi_N(ut)}{b_N-\varphi_N(ut)\varphi_{N+1}(ut)}}\)^{Z^+_N}}\).\nonumber
\end{eqnarray}

\smallskip
\noindent
Replacing $\varphi_N\varphi_{N+1}$ of the first line by the quantity $\displaystyle\frac{a_N\varphi_{N-1}\varphi_N}{b_N-\varphi_N\varphi_{N+1}}$, we iterate and obtain~:

$$\chi_D^{+,N}=\frac{\varphi_1a_1}{b_1-\frac{\varphi_1\varphi_2a_2}{b_2\cdots-\frac{\varphi_{N-1}\varphi_Na_N}{b_N-\varphi_N\varphi_{N+1}}}}.$$

\medskip
\noindent
Dividing by $\varphi_1,\cdots,\varphi_N$ at each successive level, using that the $\varphi_n$ are close to $1$, hence not $0$, uniformly in $n$ and $u\in S_+^{d-1}$ for small $t$, we get~:

$$\chi_D^{+,N}(ut)=[(a_1,b_1/\varphi_1(ut));(-a_2,b_2/\varphi_2(ut));\cdots;(-a_N,b_N/\varphi_N(ut)-\varphi_{N+1}(ut))].$$

\smallskip
\noindent
Now in \eqref{chiinit}, $\chi_D^{+,N}$ converges pointwise to $\chi_D$ by dominated convergence. Hence, by lemma \eqref{cf}~:

$$\chi_D^{+}(ut)=[(a_1,b_1/\varphi_1(ut));(-a_2,b_2/\varphi_2(ut));\cdots;(-a_n,b_n/\varphi_n(ut))\cdots].$$

\smallskip
\noindent
A similar expression is true for $\chi_D^-(ut)$. We have in fact shown something slightly stronger~:
\begin{lemme}

$ $

\label{plusgen}

\noindent
Let $(\gamma_n)_{n\geq1}$ be a sequence of complex numbers with $0<|\gamma_n|\leq1$. Then~:

$$\E^{+}\prod_{k=1}^{\sigma-1}\gamma_{Y_k}=\gamma_1\E^+\prod_{n\geq 1}[\gamma_n\gamma_{n+1}]^{Z^+_{n+1}}=[(a_1,b_1/\gamma_1);(-a_2,b_2/\gamma_2);\cdots;(-a_n,b_n/\gamma_n);\cdots].$$
\end{lemme}

\subsection{Another reduction}

Let $E_n(ut)=\sum_{k\in\Z^d}e^{itu.k}\mu_n(k)$, $t\in\R$. From \eqref{varphin}, $\varphi_n(ut)=(1-r_n)/(1-r_nE_n(ut))$, giving~:

$$\frac{1}{\varphi_n(ut)}=1-itu.m_n\frac{r_n}{1-r_n}+O(t^2),$$

\smallskip
\noindent
with $O$ uniform in $n$ and $u\in S_+^{d-1}$. We shall replace below the $\varphi_n(ut)$ by the $\psi_n(ut)$ in the recursive relation \eqref{syst} satisfied by the $(B_n)$, where~:

\begin{defi}

$ $

\noindent
For $n\in\Z$, $u\in S_+^{d-1}$ and $t\in\R$, set $\eta_n=r_nm_n/p_n$ and $\displaystyle\frac{1}{\psi_n(ut)}=1-itu.\eta_n\frac{p_n}{1-r_n}=1-itu.\eta_n/b_n$.

\end{defi}

\begin{lemme}
$  $

\noindent
Let $c=\delta^3/4>0$. For small $t>0$, uniformly in $n$ and $u\in S_+^{d-1}$~:

\begin{equation}
\label{minorphi}
|\varphi_n(ut)|\leq 1-ct^2.
\end{equation} 

\end{lemme}

\smallskip
\noindent
\textit{Proof of the lemma~:} 

\noindent
Let $M_{2,n}(u)=\sum_{k\in\Z^d}(k.u)^2\mu_n(k)$, $m_n(u)=m_n.u$ and $Var_n(u)=M_{2,n}(u)-m_n(u)^2$. A computation gives~:

$$|\varphi_n(ut)|=1-\frac{t^2}{2}\frac{r_n}{(1-r_n)^2}(M_{2,n}(u)-r_nVar_n(u))+O(t^3),$$

\noindent
with $O$ uniform in $n$ and $u\in S^{d-1}$, due to the uniformly bounded third moment of $\mu_n$. Using the hypotheses, we have $\delta^2\leq\delta M_{2,n}(u)\leq M_{2,n}(u)-r_nVar_n(u).$ Hence~:

$$|\varphi_n(ut)|\leq1-\frac{t^2\delta^3}{2}+O(t^3)\leq1-\frac{t^2\delta^3}{4},$$ 

\noindent
for $t$ small enough, uniformly in $n$ and $u\in S_+^{d-1}$.\fin

\begin{lemme}
\label{error}

$ $

\noindent
Let $R^+(t)=1-\E^+((1-t^2)^{\sigma-1})$ and $f^+(ut)=\E^+(\prod_{k=1}^{\sigma-1}\psi_{Y_k}(ut))$. 

\begin{enumerate}
\item For all $C\geq 1$, for $x>0$ large enough~: $\psi_+^{-1}(Cx)\leq 2C^2\psi_+^{-1}(x)$. 

\item There exists $\alpha\geq1$ so that for small $t>0$~:

$$\frac{1}{\alpha}\leq R^+(t)\psi_+^{-1}(1/t)\leq\alpha.$$

\item 
There exist constants $C_1>0,C_2>0$ so that for small $t>0$, uniformly in $u\in S_+^{d-1}$~:

$$1-|\chi_D^+(ut)|\geq C_1R^+(t)\mbox{ and }\left|{\chi_D^+(ut)-f^+(ut)}\right|\leq C_2R^+(t).$$

 \end{enumerate}
\end{lemme}

\smallskip
\noindent
\textit{Proof of the lemma~:} 

\noindent
$1.$ Recall that $\psi_+^2(n)=nw_+\circ v_+^{-1}(n)$, so $n\longmapsto \psi^2_+(n)/n$ is non-decreasing. Let $C\geq 1$ and $x>0$. Set $n=\psi_+^{-1}(x)$ and suppose that $n\geq1$. By definition, $\psi_+(n)\leq x<\psi_+(n+1)$. Similarly, let $n+p=\psi_+^{-1}(Cx)$. Then~:

$$\frac{\psi_+^{-1}(Cx)}{\psi_+^{-1}(x)}=\frac{n+p}{n}\leq2\frac{n+p}{n+1}\leq2\frac{\psi_+^2(n+p)}{\psi_+^2(n+1)}\leq 2\frac{C^2x^2}{x^2}=2C^2.$$

\medskip
\noindent
$2.$ As a preliminary point, for $n\geq1$, let $\Theta_+(n)>0$ be such that $\Theta_+^2(n)=\sum_{1\leq k\leq l\leq n}(\rho_l/\rho_k)$. Fix constants $c>0,c'>0$ so that $c w_+(k)\leq\sum_{1\leq u\leq k}(1/\rho_u)$ and $v_+(k+1)\leq c'v_+(k)$, for all $k\geq1$. We claim that there exists $C>0$ so that for all $x>0$ large enough~:

$$(1/C)v_+\circ\Theta_+^{-1}(x)\leq \psi_+^{-1}(x)\leq Cv_+\circ\Theta_+^{-1}(x).$$

\smallskip
\noindent
The second inequality follows from the remark that $\Theta_+\circ v_+^{-1}\leq \psi_+$, giving $v_+^{-1}\circ\psi_+^{-1}\leq \Theta_+^{-1}$, and the fact that $v_+(v_+^{-1}(x))\geq c''x$, for some constant $c''>0$. For the first one, let $x>0$ and $n=\Theta_+^{-1}(x)$. For any $1\leq m\leq n$~:

$$x^2\geq (v_+(n)-v_+(m))\sum_{1\leq k\leq m}(1/\rho_k)\geq c(v_+(n)-v_+(m))w_+(m).$$

\smallskip
\noindent
Choose $m\leq n$ so that $v_+(m)\leq v_+(n)/2<v_+(m+1)$. Hence, as $m=v_+^{-1}(v_+(m))$~:

$$x^2\geq(c/2)v_+(m)w_+(m)\geq(c/2)\psi_+^2(v_+(m)).$$

\smallskip
\noindent
We obtain, using at the end the first point of the lemma~:

$$v_+(n)/(2c')\leq v_+(m)\leq \psi_+^{-1}(x/\sqrt{c/2})\leq (4/c)\psi_+^{-1}(x).$$

\smallskip
\noindent
This completes the proof of the claim. 

\medskip
Let us now turn to the evaluation of $R^+(t)$. Using lemmas \ref{plusgen} and \ref{cf} we have~:

$$\E^+((1-t^2)^{\sigma-1})=\lim_{n\rightarrow+\infty}\frac{\alpha_n(t)}{\beta_n(t)}=\sum_{n\geq 1}\frac{\rho_n}{\beta_n\beta_{n-1}},$$

\noindent
where $\beta_{-1}=0$, $\beta_0=1$ and $\beta_n=(b_n/(1-t^2))\beta_{n-1}-a_n\beta_{n-2}$. We omit the dependence in $t$. The $(\alpha_n)$ satisfy the same recursive relation with this time $\alpha_{-1}=1$ and $\alpha_0=0$. First, as there is a constant $C>0$ so that for all $n\geq1$, $\Theta_+(n)\leq\Theta_+(n+1)\leq C\Theta_+(n)$, we deduce that for any constant $c>0$ (chosen later), there exists a constant $c'>0$ so that for small enough $t>0$ there is an integer $N(t)$ so that~:

$$\frac{c'}{t^2}\leq \Theta^2_+(N(t))\leq \frac{c}{t^2}.$$

\noindent
Next, using lemma \ref{cf}~:

\begin{eqnarray}
\label{troistermes}
\left|{R^+(t)- \({\frac{\beta_{N(t)}(t)-\alpha_{N(t)}(t)}{\beta_{N(t)}(t)}}\)}\right|\leq\sum_{n>N(t)}\frac{\rho_n}{\beta_n(t)\beta_{n-1}(t)}\leq \frac{1}{\beta_{N(t)}}.
\end{eqnarray}

\noindent
We shall show that there exists a constant $\varepsilon>0$ so that $1+\varepsilon\leq \beta_{N(t)}(t)-\alpha_{N(t)}(t)\leq 1/\varepsilon$ and next that $v_+(N(t))\leq\beta_{N(t)}(t)\leq v_+(N(t))/\varepsilon$. These two properties imply that $R^+(t)$ has exact order $1/v_+(N(t))$ and so $1/\psi_+^{-1}(1/t)$, by the claim and the first point. 

\medskip
\noindent
We have $b_n/(1-t^2)=b_n+t^2c_n(t)$, with $(1/\alpha)\leq c_n(t)\leq \alpha$, for some constant $\alpha>0$. Next~:

$$
\begin{pmatrix}
\beta_n\\
\beta_{n-1}\\
\end{pmatrix}=\begin{pmatrix}
b_n+t^2c_n(t)&-a_n\\
1&0\\
\end{pmatrix}\cdots\begin{pmatrix}
b_1+t^2c_1(t)&-a_1\\
1&0\\
\end{pmatrix} \begin{pmatrix}
1\\
0\\
\end{pmatrix}.$$

\medskip
\noindent
Setting $~
C_n=\begin{pmatrix}
b_n&-a_n\\
1&0\\
\end{pmatrix}$, $B=\begin{pmatrix}
1&0\\
0&0\\
\end{pmatrix}$ and since $\beta_n(0)=v_+(n)$, we obtain~:

\begin{eqnarray}
\beta_n&=&v_+(n)+\sum_{r=1}^nt^{2r}\sum_{1\leq k_1<\cdots<k_r\leq n}c_{k_1}(t)\cdots c_{k_r}(t)\langle e_1,C_n\cdots C_{k_r+1}B\cdots BC_{k_1-1}\cdots C_1 e_1\rangle\nonumber\\
&=&v_+(n)+\sum_{r=1}^nt^{2r}\sum_{1\leq k_1<\cdots<k_r\leq n}(c_{k_1}\cdots c_{k_r})(t)v_+(k_1-1)\theta^{k_1}v_+(k_2-k_1-1)\cdots \theta^{k_r}v_+(n-k_r).\nonumber
\end{eqnarray}

\noindent
Idem, since $\alpha_n=a_1\theta\beta_{n-1}$~:

\begin{eqnarray}
\alpha_n&=&v_+(n)-1+\sum_{r=1}^nt^{2r}\sum_{2\leq k_1<\cdots<k_r\leq n}(c_{k_1}\cdots c_{k_r})(t)(v_+(k_1-1)-1)\theta^{k_1}v_+(k_2-k_1-1)\cdots \theta^{k_r}v_+(n-k_r).\nonumber
\end{eqnarray}

\noindent
This furnishes~:

\begin{eqnarray}\beta_n-\alpha_n&=&1+\sum_{r=1}^nt^{2r}\sum_{1\leq k_1<\cdots<k_r\leq n}(c_{k_1}\cdots c_{k_r})(t)\theta^{k_1}v_+(k_2-k_1-1)\cdots \theta^{k_r}v_+(n-k_r).\nonumber
\end{eqnarray}

\noindent
As a result $\beta_n\leq v_+(n)(1+\sum_{1\leq r\leq n}\alpha^rt^{2r}(\Theta^2_+(n))^r)$ and $\beta_n-\alpha_n\geq1+t^{2}\Theta^2_+(n)/\alpha$. We simply choose $0<c\leq \alpha/2$ to get the desired result.

\medskip
\noindent
$3.$ We have $\chi_D^+(ut)=\E^+(\prod_{k=1}^{\sigma-1}\varphi_{Y_k}(ut))$. By \eqref{minorphi}, $|\chi_D^+(ut)|\leq \E^+((1-ct^2)^{\sigma-1})$. This gives the first inequality, as the first point of the lemma says that $R^+(\sqrt{c}t)\leq CR^+(t)$, for some constant $C$ depending on $c$. Concerning the second inequality~:

\begin{eqnarray}
\left|{\chi_D^+(ut)-f^+(ut)}\right|&=&\left|{\E^+\({\prod_{k=1}^{\sigma-1}\varphi_{Y_k}(ut)}\)-\E^+\({\prod_{k=1}^{\sigma-1}\psi_{Y_k}(ut)}\)}\right|\nonumber\\
&=&\left|{\E^+\({\sum_{k=1}^{\sigma-1}\({\prod_{l=1}^{k-1}\varphi_{Y_l}(ut)(\varphi_{Y_k}(ut)-\psi_{Y_k}(ut))\prod_{l=k+1}^{\sigma-1}\psi_{Y_l}(ut)}\)}\)}\right|.\nonumber\\
&\leq&\E^+\({\sum_{k=1}^{\sigma-1}\({\prod_{l=1}^{k-1}|\varphi_{Y_l}(ut)||\varphi_{Y_k}(ut)-\psi_{Y_k}(ut)|\prod_{l=k+1}^{\sigma-1}|\psi_{Y_l}(ut)|}\)}\).\nonumber\end{eqnarray}

\noindent
Using now that for some $C>0$ and small enough $t>0$, uniformly in $n$ and $u\in S_+^{d-1}$, $|\varphi_n(ut)-\psi_n(ut)|\leq Ct^2$, as well as $|\varphi_n(ut)|\leq 1-ct^2$ and $|\psi_n(ut)|\leq1$, we get for small $t>0$~:

$$\left|{\chi_D^+(ut)-f^+(ut)}\right|\leq Ct^2\E^+\({\sum_{k=1}^{\sigma-1}(1-ct^2)^{k-1}}\)=Ct^2\E^+\({\frac{1-(1-ct^2)^{\sigma-1}}{ct^2}}\)=\frac{C}{c}R^+(\sqrt{c}t).$$

\noindent
The conclusion now comes from the first point of the lemma. 

\fin
\section{Precise analysis of some convergents}
As a summary, from the previous section, uniformly in $u\in S_+^{d-1}$~:

$$\chi^+_D(ut)=f^+(ut)+O(R^+(t)),$$

\smallskip
\noindent
with $f^+(ut)=\lim_{n\rightarrow+\infty}A_n(ut)/B_n(ut)$, where now~:

$$
\begin{pmatrix}
B_n(ut)\\
B_{n-1}(ut)\\
\end{pmatrix}=\begin{pmatrix}
b_n-itu.\eta_n&-a_n\\
1&0\\
\end{pmatrix}\cdots\begin{pmatrix}
b_1-itu.\eta_1&-a_1\\
1&0\\
\end{pmatrix} \begin{pmatrix}
1\\
0\\
\end{pmatrix},$$

\medskip
\noindent
together with $A_n(ut)=a_1\theta B_{n-1}(ut)$.

\medskip
Recall the definitions $R_k^l(u)=\sum_{k\leq r\leq l}\eta_r.u(\rho_l/\rho_r)$ and $T_k^l(u)=(R_k^l(u))^2\rho_{k-1}/\rho_l$, $k\leq l$. For fixed $u\in S_+^{d-1}$, notice that these quantities depend only on the data in $[k,l]$.

\begin{defi}
$ $

\noindent 
We fix $u\in S_+^{d-1}$. Omitting the dependence with respect to $u$, set~:

$$\Delta_r^n=\sum_{1\leq k_1<\cdots<k_r\leq n}R_1^{k_1}R_{k_1+1}^{k_2}\cdots R_{k_{r-1}+1}^{k_r},$$

\noindent
with $\Delta_0^n=1$ and $\Delta_r^n=0$ if $r>n$ or $r<0$.
\end{defi}

\smallskip
Proceeding as in the previous section, setting $\eta_k'=u.\eta_k$, we develop~:

\begin{eqnarray}
B_n(ut)&=&v_+(n)+\sum_{r=1}^n(-it)^{r}\sum_{1\leq k_1<\cdots<k_r\leq n}\eta'_{k_1}\cdots \eta'_{k_r}v_+(k_1-1)\theta^{k_1}v_+(k_2-k_1-1)\cdots \theta^{k_r}v_+(n-k_r),\nonumber\\
A_n(ut)&=&v_+(n)-1+\sum_{r=1}^n(-it)^{r}\sum_{2\leq k_1<\cdots<k_r\leq n}\eta'_{k_1}\cdots \eta'_{k_r}(v_+(k_1-1)-1)\theta^{k_1}v_+(k_2-k_1-1)\cdots \theta^{k_r}v_+(n-k_r).\nonumber
\end{eqnarray}

\noindent
We therefore obtain~:

$$B_n(ut)-A_n(ut)=1+\sum_{r=1}^n(-it)^{r}\sum_{1\leq k_1<\cdots<k_r\leq n}\eta'_{k_1}\cdots \eta'_{k_r}\theta^{k_1}v_+(k_2-k_1-1)\cdots \theta^{k_r}v_+(n-k_r).$$

\medskip
\noindent
In the last sum, fix $k_2,\cdots,k_r$ and write~:

$$\sum_{1\leq k_1<k_2}\eta'_{k_1}\theta^{k_1}v_+(k_2-k_1-1)=\sum_{1\leq k_1<k_2}\eta'_{k_1}\sum_{k_1\leq l<k_2}\frac{\rho_l}{\rho_{k_1}}=\sum_{1\leq l<k_2}\sum_{1\leq k_1\leq l}\eta'_{k_1}\frac{\rho_l}{\rho_{k_1}}=\sum_{1\leq l<k_2}R_{1}^l.$$

\noindent
Successively iterate this manipulation for $k_2,\cdots,k_r$ in the formula for $B_n(ut)-A_n(ut)$. Then~:

$$B_n(ut)-A_n(ut)=1+\sum_{r=1}^n(-it)^{r}\sum_{1\leq k_1<\cdots<k_r\leq n}R_1^{k_1}R_{k_1+1}^{k_2}\cdots R_{k_{r-1}+1}^{k_r}=\sum_{r=0}^n(-it)^{r}\Delta_r^n.$$

\noindent
Similarly, using as first step that $\sum_{1\leq k_1<k_2}\eta'_{k_1}v_+(k_1-1)\theta^{k_1}v_+(k_2-k_1-1)=\sum_{0\leq s<l<k_2}\rho_sR_{s+1}^l$~:

$$B_n(ut)=v_+(n)+\sum_{r=1}^n(-it)^{r}\sum_{0\leq k_1<\cdots<k_{r+1}\leq n}\rho_{k_1}R_{k_1+1}^{k_2}\cdots R_{k_{r}+1}^{k_{r+1}}=\sum_{r=0}^n(-it)^{r}\sum_{0\leq k\leq n}\rho_k\theta^k\Delta_r^{n-k}.$$

\begin{prop}
\label{estim}

$ $

\noindent
Set $2^{(k,l)}=2$ if $k\not=l$ and $1$ if $k=l$. We have the following exact computations~:

\begin{enumerate}
\item $|B_n(ut)-A_n(ut)|^2=\sum_{r=0}^nt^{2r}K_r(n)$, with~:

$$K_r(n)=\sum_{1\leq l_1<k_2\leq l_2<\cdots<k_r\leq l_r< k_{r+1}\leq n+1}T_1^{l_1}\cdots T_{k_{r}}^{l_{r}}\rho_{k_{r+1}-1}2^{H_r((k_i),(l_j))},$$

where $H_r((k_i),(l_j)):=\#\{1\leq i\leq r~|~l_i+1<k_{i+1}\}$.

\item $|B_n(ut)|^2=\sum_{r=0}^nt^{2r}L_r(n)$, with $L_r(n)=\sum_{0\leq k\leq l\leq n}\rho_k\rho_l2^{(k,l)}\theta^lK_r(n-l)$.

\item $\mbox{Re}((B_n-A_n)\bar{B}_n)(ut)=\sum_{r=0}^nt^{2r}M_r(n)$, with $M_r(n)=\sum_{0\leq k\leq n}\rho_k\theta^kK_r(n-k)$.

\item $\mbox{Im}(A_n(ut)\bar{B}_n(ut))=\sum_{r=0}^{n-1}t^{2r+1}N_r(n)$, with $N_r(n)=\sum_{1\leq k\leq l\leq n}R_1^k2^{(k,l)}\rho_l\theta^lK_r(n-l)$.

\end{enumerate}
\noindent
When $r>n$ or $r<0$, set $K_r(n)=L_r(n)=M_r(n)=0$. Idem $N_r(n)=0$, $r\geq n$ or $r<0$. 

\end{prop}

\noindent
\begin{remark}
Recall that $R_k^l$ and $T_k^l$ and therefore $K_r(n)$, $L_r(n)$, $M_r(n)$, $N_r(n)$ depend on $u\in S_+^{d-1}$ but that the dependence is omitted in the notations.
\end{remark}

\medskip
\noindent
\textit{Proof of the proposition~:}

\noindent
$1.$ Since $B_n(ut)-A_n(ut)=\sum_{0\leq r\leq n}(-it)^r\Delta_r^n$, this gives~:

$$|B_n(ut)-A_n(ut)|^2=(B_n(ut)-A_n(ut))\overline{(B_n(ut)-A_n(ut))}
=\sum_{r=0}^nt^{2r}\sum_{p=-r}^r\Delta_{r+p}^n\Delta_{r-p}^n(-i)^{r+p}i^{r-p},$$

\noindent
using the conventions for $\Delta_r^n$ concerning the value of $r$ with respect to $n$. Hence $|B_n(ut)-A_n(ut)|^2=\sum_{r=0}^nt^{2r}K_r(n)$, with $K_0(n)=1$ and~:

$$K_r(n)=\sum_{p=-r}^r(-1)^p\Delta_{r+p}^n\Delta_{r-p}^n,~r\geq 1.$$

\smallskip
\noindent
We will show that~:

\begin{equation}
\label{init}
K_1(n)=\sum_{1\leq k\leq l\leq n}T_1^k\rho_l2^{(k,l)},
\end{equation}
 
\smallskip
\noindent 
together with the following recursive relation, for $r\geq 2$~:

\begin{equation}
\label{recur}
K_r(n)=\sum_{1\leq k\leq l\leq n}T_1^k\rho_l\theta^{l}K_{r-1}(n-l)2^{(k,l)}.
\end{equation}

\medskip
This then gives the announced formula. For the initial relation~:

\begin{eqnarray}K_1(n)&=&(\Delta_1^n)^2-2\Delta_2^n=\({\sum_{1\leq k\leq n}R_1^k}\)^2-2\sum_{1\leq k<l\leq n}R_1^kR_{k+1}^l\nonumber\\
&=&\sum_{1\leq k\leq n}(R_1^k)^2+2\sum_{1\leq k<l\leq n}R_1^k(R_1^l-R_{k+1}^l).\nonumber\end{eqnarray}

\noindent
Observing that $R_1^k(R_1^l-R_{k+1}^l)=(R_1^k)^2(\rho_l/\rho_k)=T_1^k\rho_l$, this proves \eqref{init}. Let us now turn to the proof of \eqref{recur}. Taking first general $p\geq1$ and $q\geq 1$, we write~:

$$\Delta_p^n\Delta_q^n=\sum_{\substack{1\leq k_1<\cdots<k_p\leq n\\1\leq k_1'<\cdots<k_q'\leq n}}(R_1^{k_1}\cdots R_{k_{p-1}+1}^{k_p})(R_1^{k'_1}\cdots R_{k'_{q-1}+1}^{k'_q}).$$

\noindent
Distinguishing the cases $k_1=k_1'$, $k_1<k_1'$ and $k'_1<k_1$, we decompose~:

\begin{eqnarray}
\Delta_p^n\Delta_q^n=&~&\sum_{1\leq k\leq n}(R_1^{k})^2\theta^{k}\Delta_{p-1}^{n-k}\theta^{k}\Delta_{q-1}^{n-k}\nonumber\\
&+&\sum_{\substack{1\leq k_1<\cdots<k_p\leq n\\k_1< k_1'<\cdots<k_q'\leq n}}R_1^{k_1}(R_{k_1+1}^{k_2}\cdots R_{k_{p-1}+1}^{k_p})(R_1^{k_1}\frac{\rho_{k'_1}}{\rho_{k_1}}+R_{k_1+1}^{k_1'})(R_{k_1'+1}^{k_2'}\cdots R_{k'_{q-1}+1}^{k'_q})\nonumber\\
&+&\sum_{\substack{1\leq k'_1<\cdots<k'_q\leq n\\k'_1< k_1<\cdots<k_p\leq n}}(R_1^{k'_1}\frac{\rho_{k_1}}{\rho_{k'_1}}+R_{k_1'+1}^{k_1})(R_{k_1+1}^{k_2}\cdots R_{k_{p-1}+1}^{k_p})(R_1^{k'_1}\cdots R_{k'_{q-1}+1}^{k'_q}).\nonumber\end{eqnarray}

\noindent
Regrouping terms, this is rewritten as~:

\begin{eqnarray}
\Delta_p^n\Delta_q^n&=&\sum_{1\leq k\leq n}(R_1^{k})^2\[{\theta^k\Delta_{p-1}^{n-k}\sum_{k\leq l\leq n}\theta^l\Delta_{q-1}^{n-l}\frac{\rho_l}{\rho_k}+\theta^k\Delta_{q-1}^{n-k}\sum_{k<l\leq n}\theta^l\Delta_{p-1}^{n-l}\frac{\rho_l}{\rho_k}}\]\nonumber\\
&+&\sum_{1\leq k\leq n}R_1^k\[{\theta^k\Delta_{p-1}^{n-k}\theta^k\Delta_q^{n-k}+\theta^k\Delta_{p}^{n-k}\theta^k\Delta_{q-1}^{n-k}}\].\nonumber
\end{eqnarray}

\noindent
Taking $r\geq 2$, insert the latter in $K_r(n)=\sum_{-r+1\leq p\leq r-1}(-1)^p\Delta_{r+p}^n\Delta_{r-p}^n+2(-1)^r\Delta_{2r}^n$ and get~:

\begin{eqnarray}
K_r(n)&=&\sum_{1\leq k\leq n}(R_1^k)^2\sum_{-r+1\leq p\leq r-1}(-1)^p\[{\theta^k\Delta_{r+p-1}^{n-k}\sum_{k\leq l\leq n}\theta^l\Delta_{r-p-1}^{n-l}\frac{\rho_l}{\rho_k}+\theta^k\Delta_{r-p-1}^{n-k}\sum_{k<l\leq n}\theta^l\Delta_{r+p-1}^{n-l}\frac{\rho_l}{\rho_k}}\]\nonumber\\
&+&2(-1)^r\Delta_{2r}^n+2\sum_{1\leq k\leq n}R_1^k\[{\sum_{-r+1\leq p\leq r-1}(-1)^{p}\theta^k\Delta_{r+p-1}^{n-k}\theta^k\Delta_{r-p}^{n-k}}\].\nonumber
\end{eqnarray}

\noindent
The last line is $2\sum_{1\leq k\leq n}R_1^k[\sum_{-r+1\leq p\leq r}(-1)^{p}\theta^k\Delta_{r+p-1}^{n-k}\theta^k\Delta_{r-p}^{n-k}]$. The bracketed sum is $0$, for instance when doing the change of variable $p\longmapsto -p+1$. Separating now the term with $k=l$ in the first sum above and recognizing $\theta^kK_{r-1}(n-k)$, we obtain~:

\begin{eqnarray}
K_r(n)&=&\sum_{1\leq k\leq n}(R_1^k)^2\[{\theta^kK_{r-1}(n-k)+2\sum_{-r+1\leq p\leq r-1}(-1)^p\theta^k\Delta_{r+p-1}^{n-k}\sum_{k<l\leq n}\theta^l\Delta_{r-p-1}^{n-l}\frac{\rho_l}{\rho_k}}\].\nonumber\end{eqnarray}

\noindent
Setting $m=n-k$ and $Z_r(m)=\sum_{-r\leq p\leq r}(-1)^p\Delta_{r+p}^m\sum_{1\leq l\leq m}\theta^l\Delta_{r-p}^{m-l}\rho_l$, we therefore have~:

$$K_r(n)=\sum_{1\leq k\leq n}(R_1^k)^2\[{\theta^kK_{r-1}(n-k)+2\theta^kZ_{r-1}(n-k)}\].$$

\smallskip
We shall show that~:

\begin{equation}
\label{zmfinal}
Z_r(m)=\sum_{1\leq k\leq m}\theta^kK_{r}(m-k)\rho_k,~r\geq 1.
\end{equation}

\noindent
To complete the proof of (\ref{recur}), we simply apply this to $Z_{r-1}(n-k)$ in the previous equality. First of all, with $0\leq p\leq r-1$~:

\begin{eqnarray}
\Delta_{r+p}^m\sum_{1\leq l\leq m}\theta^l\Delta_{r-p}^{m-l}\rho_l&=&\sum_{\substack{1\leq k_1<\cdots<k_{r+p}\leq m\\1\leq l_1<l_2<\cdots<l_{r-p+1}\leq m}}R_1^{k_1}\cdots R_{k_{r+p-1}+1}^{k_{r+p}}R_{l_1+1}^{l_2}\cdots R_{l_{r-p}+1}^{l_{r-p+1}}\rho_{l_1}\nonumber\\
&=&\sum_{\substack{1\leq k_1<\cdots <k_{r+p}\leq m\\k_1\leq l_1<\cdots<l_{r-p+1}\leq m}}R_1^{k_1}R_{k_1+1}^{k_2}\cdots R_{k_{r+p-1}+1}^{k_{r+p}}R_{l_1+1}^{l_2}\cdots R_{l_{r-p}+1}^{l_{r-p+1}}\rho_{l_1}\nonumber\\
&+&\sum_{\substack{1\leq l_1<\cdots <l_{r-p+1}\leq m\\l_1<k_1<\cdots<k_{r+p}\leq m}}(R_1^{l_1}\frac{\rho_{k_1}}{\rho_{l_1}}+R_{l_1+1}^{k_1})R_{k_1+1}^{k_2}\cdots R_{k_{r+p-1}+1}^{k_{r+p}}R_{l_1+1}^{l_2}\cdots R_{l_{r-p}+1}^{l_{r-p+1}}\rho_{l_1}.\nonumber
\end{eqnarray}

\noindent
Written in a more concise way~:

\begin{eqnarray}
\Delta_{r+p}^m\sum_{1\leq l\leq m}\theta^l\Delta_{r-p}^{m-l}\rho_l&=&\sum_{1\leq k\leq m}R_1^k\[{\theta^k\Delta_{r+p-1}^{m-k}\sum_{k\leq l\leq m}\theta^l\Delta_{r-p}^{m-l}\rho_l+\theta^k\Delta_{r-p}^{m-k}\sum_{k<l\leq m}\theta^l\Delta_{r+p-1}^{m-l}\rho_l}\]\nonumber\\
&+&\sum_{1\leq k\leq m}\theta^k\Delta_{r+p}^{m-k}\theta^k\Delta_{r-p}^{m-k}\rho_k.\nonumber
\end{eqnarray}

\noindent
This allows to write~:

\begin{eqnarray}
Z_r(m)&=&(-1)^{r}\[{\Delta_{2r}^m\sum_{1\leq l\leq m}\rho_l+\sum_{1\leq l\leq m}\theta^l\Delta_{2r}^{m-l}\rho_l}\]+\sum_{1\leq k\leq m}\sum_{-r+1\leq p\leq r-1}(-1)^p\theta^k\Delta_{r+p}^{m-k}\theta^k\Delta_{r-p}^{m-k}\rho_k\nonumber\\
&+&\sum_{1\leq k\leq m}R_1^k\sum_{-r+1\leq p\leq r-1}(-1)^p\[{\theta^k\Delta_{r+p-1}^{m-k}\sum_{k\leq l\leq m}\theta^l\Delta_{r-p}^{m-l}\rho_l+\theta^k\Delta_{r-p}^{m-k}\sum_{k<l\leq m}\theta^l\Delta_{r+p-1}^{m-l}\rho_l}\].\nonumber
\end{eqnarray}

\noindent
Recognizing some $\theta^kK_r(m-k)$, we get~:

\begin{eqnarray}
Z_r(m)&=&(-1)^{r}\[{\Delta_{2r}^m\sum_{1\leq l\leq m}\rho_l-\sum_{1\leq l\leq m}\theta^l\Delta_{2r}^{m-l}\rho_l}\]+\sum_{1\leq k\leq m}\theta^kK_r(m-k)\rho_k\nonumber\\
&+&\sum_{1\leq k\leq m}R_1^k\sum_{-r+1\leq p\leq r-1}(-1)^p\[{\theta^k\Delta_{r+p-1}^{m-k}\sum_{k\leq l\leq m}\theta^l\Delta_{r-p}^{m-l}\rho_l}\]\nonumber\\
&+&\sum_{1\leq k\leq m}R_1^k\sum_{-r+2\leq p\leq r}(-1)^{p+1}\[{\theta^k\Delta_{r+p-1}^{m-k}\sum_{k< l\leq m}\theta^l\Delta_{r-p}^{m-l}\rho_l}\].\nonumber\end{eqnarray}

\noindent
Consequently~:

\begin{eqnarray}
Z_r(m)&=&\sum_{1\leq k\leq m}\theta^kK_r(m-k)\rho_k+(-1)^{r}\[{\Delta_{2r}^m\sum_{1\leq l\leq m}\rho_l-\sum_{1\leq l\leq m}\theta^l\Delta_{2r}^{m-l}\rho_l}\]\nonumber\\
&+&(-1)^{r+1}\sum_{1\leq k\leq m}R_1^k\({\theta^k\Delta_{2r-1}^{m-k}\sum_{k\leq l\leq m}\rho_l+\sum_{k<l\leq m}\theta^l\Delta_{2r-1}^{m-l}\rho_l}\)\nonumber\\
&+&\sum_{1\leq k\leq m}R_1^k\sum_{-r+1\leq p\leq r}(-1)^p\[{\theta^k\Delta_{r+p-1}^{m-k}\sum_{k\leq l\leq m}\theta^l\Delta_{r-p}^{m-l}\rho_l-\theta^k\Delta_{r+p-1}^{m-k}\sum_{k< l\leq m}\theta^l\Delta_{r-p}^{m-l}\rho_l}\].\nonumber\end{eqnarray}

\noindent
The last line is $\sum_{1\leq k\leq m}R_1^k[\sum_{-r+1\leq p\leq r}(-1)^p\theta^k\Delta_{r+p-1}^{m-k}\theta^k\Delta_{r-p}^{m-k}\rho_k]$. For the same reason as before, the inside brackets are $0$. Therefore it finally remains to show that the sum of the second and third terms is also $0$, in other words that~:

$$\Delta_{2r}^m\sum_{1\leq l\leq m}\rho_l-\sum_{1\leq l\leq m}\theta^l\Delta_{2r}^{m-l}\rho_l-\sum_{1\leq k\leq m}R_1^k\({\theta^k\Delta_{2r-1}^{m-k}\sum_{k\leq l\leq m}\rho_l+\sum_{k<l\leq m}\theta^l\Delta_{2r-1}^{m-l}\rho_l}\)=0.$$

\noindent
Equivalently~:
$$\sum_{1\leq k\leq m}R_1^k\theta^k\Delta_{2r-1}^{m-k}\sum_{1\leq k<l}\rho_l-\sum_{1\leq l\leq m}\theta^l\Delta_{2r}^{m-l}\rho_l-\sum_{1\leq k\leq m}R_1^k\sum_{k<l\leq m}\theta^l\Delta_{2r-1}^{m-l}\rho_l=0.$$

\noindent
In the last term, replace $R_1^k$ by $(R_1^l-R_{k+1}^l)\rho_k/\rho_l$. It remains to show that~:

$$-\sum_{1\leq l\leq m}\theta^l\Delta_{2r}^{m-l}\rho_l+\sum_{1\leq k\leq m}R_1^k\theta^k\Delta_{2r-1}^{m-k}\sum_{1\leq l<k}\rho_l-\sum_{1\leq k<l\leq m}R_1^l\theta^l\Delta_{2r-1}^{m-l}\rho_k+\sum_{1\leq k<l\leq m}R_{k+1}^l\theta^l\Delta_{2r-1}^{m-l}\rho_k=0.$$

\noindent
As this is true, this completes the proof of this first point.

\medskip
\noindent
$2.$ Let us define $\tilde{\Delta}_r^n=\sum_{0\leq k\leq n}\rho_k\theta^k\Delta_r^{n-k}$, so that $B_{n}(ut)=\sum_{0\leq r\leq n}(-it)^r\tilde{\Delta}_r^n$. As for $|B_n(ut)-A_n(ut)|^2$ in the first point, we have~:

$$|B_n(ut)|^2=\sum_{0\leq r\leq n}t^{2r}L_r(n),~\mbox{where }L_r(n)=\sum_{-r\leq p\leq r}(-1)^{p}\tilde{\Delta}_{r+p}^n\tilde{\Delta}_{r-p}^n.$$

\noindent
In order to compute $L_r(n)$, notice first that~:

\begin{eqnarray}
\tilde{\Delta}_{r+p}^n\tilde{\Delta}_{r-p}^n&=&\sum_{0\leq k\leq n}\rho_k\[{\theta^k\Delta_{r+p}^{n-k}\sum_{k\leq l\leq n}\theta^l\Delta_{r-p}^{n-l}\rho_l+\theta^k\Delta_{r-p}^{n-k}\sum_{k<l\leq n}\theta^l\Delta_{r+p}^{n-l}\rho_l}\].\nonumber\end{eqnarray}

\noindent
Replacing in $L_r(n)$, this allows to write, using the expressions of $K_r(n)$ and $Z_r(n)$ given in (\ref{zmfinal})~:

\begin{eqnarray}
\label{eqlrn}
L_r(n)&=&\sum_{0\leq k\leq n}\rho_k\sum_{-r\leq p\leq r}(-1)^p\[{\theta^k\Delta_{r+p}^{n-k}\sum_{k\leq l\leq n}\rho_l\theta^l\Delta_{r-p}^{n-l}+\theta^k\Delta_{r-p}^{n-k}\sum_{k<l\leq n}\rho_l\theta^l\Delta_{r+p}^{n-l}}\]\nonumber\\
&=&\sum_{0\leq k\leq n}(\rho_k)^2\theta^kK_r(n-k)+2\sum_{0\leq k\leq n}(\rho_k)^2\theta^kZ_r(n-k)\\
&=&\sum_{0\leq k\leq n}\rho_k\[{\rho_k\theta^kK_r(n-k)+2\sum_{k<l\leq n}\rho_l\theta^lK_r(n-l)}\]=\sum_{0\leq k\leq n}\rho_k\theta^kK_r(n-k)\sum_{0\leq l\leq k}2^{(l,k)}\rho_l\nonumber.\end{eqnarray}

\noindent
This completes the proof of this point.

\medskip
\noindent
$3.$ Directly, we obtain~:

\begin{equation}
\label{imap}(B_n-A_n)(ut)\overline{B}_n(ut)=\sum_{0\leq r\leq n}(-it)^r\Delta_r^n\sum_{0\leq r'\leq n}(it)^{r'}\tilde{\Delta}_{r'}^n.
\end{equation}

\noindent
When developing and taking the real part, only terms with $r+r'$ even intervene. This gives~:

$$\mbox{Re}((B_n-A_n)\overline{B}_n)(ut)=\sum_{0\leq r\leq n}t^{2r}\[{\sum_{-r\leq p\leq r}(-i)^{r+p}i^{r-p}\Delta_{r+p}^n\tilde{\Delta}_{r-p}^n}\]=\sum_{0\leq r\leq n}t^{2r}M_r(n),$$

\noindent
with this time~:

$$M_r(n)=\sum_{-r\leq p\leq r}(-1)^p\Delta_{r+p}^n\tilde{\Delta}_{r-p}^n.$$

\noindent
Since $\tilde{\Delta}_r^n=\Delta_r^n+\sum_{1\leq k\leq n}\rho_k\theta^k\Delta_r^{n-k}$, using $K_r(n)$ and the value of $Z_r(n)$ in $(\ref{zmfinal})$, we have~:

\begin{eqnarray}
M_r(n)&=&K_r(n)+Z_r(n)=\sum_{0\leq k\leq n}\rho_k\theta^kK_r(n-k).\nonumber\end{eqnarray}

\noindent
This ends the proof of this point.

\medskip
\noindent
$4.$ In the same way as for $3.$, when taking the imaginary part in \eqref{imap}, only terms with $r+r'$ odd come into play. Consequently~:

$$\mbox{Im}(A_n\overline{B}_n)(ut)=-\frac{1}{i}\sum_{0\leq r\leq n-1}t^{2r+1}\[{\sum_{-r-1\leq p\leq r}(-i)^{r+p+1}i^{r-p}\Delta_{r+p+1}^n\tilde{\Delta}_{r-p}^n}\]=\sum_{0\leq r\leq n-1}t^{2r+1}N_r(n),$$

\noindent
with this time~:

$$N_r(n)=\sum_{-r-1\leq p\leq r}(-1)^{p}\Delta_{r+p+1}^n\tilde{\Delta}_{r-p}^n.$$

\noindent
Using again that $\tilde{\Delta}_r^n=\sum_{0\leq k\leq n}\rho_k\theta^k\Delta_r^{n-k}$, we get~:

$$N_r(n)=\sum_{0\leq k\leq n}\sum_{-r-1\leq p\leq r}(-1)^{p}\Delta_{r+p+1}^n\theta^k\Delta_{r-p}^{n-k}\rho_k.$$

\noindent
Notice that the term corresponding to $k=0$ equals $0$, for symmetry reasons as before. It remains~:

\begin{eqnarray}
\label{nrn}
N_r(n)&=&(-1)^{r+1}\sum_{1\leq k\leq n}\theta^k\Delta_{2r+1}^{n-k}\rho_k+\sum_{1\leq l\leq n}R_1^l\sum_{-r\leq p\leq r}(-1)^{p}\theta^l\Delta_{r+p}^{n-l}\sum_{l<k\leq n}\theta^k\Delta_{r-p}^{n-k}\rho_k\nonumber\\
&+&\sum_{1\leq k\leq n}\sum_{-r\leq p\leq r}(-1)^{p}\theta^k\Delta_{r-p}^{n-k}\rho_k\sum_{k\leq l\leq n}R_1^l\theta^l\Delta_{r+p}^{n-l}\nonumber\\
&=&\sum_{1\leq k\leq n}R_1^k\rho_k\theta^kK_r(n-k)+\sum_{1\leq l\leq n}R_1^l\rho_l\theta^lZ_r(n-l)+O_r(n),
\end{eqnarray}

\noindent
where we introduce~:

$$O_r(n)=(-1)^{r+1}\sum_{1\leq k\leq n}\theta^k\Delta_{2r+1}^{n-k}\rho_k+\sum_{1\leq k\leq n}\sum_{-r\leq p\leq r}(-1)^{p}\theta^k\Delta_{r-p}^{n-k}\rho_k\sum_{k<l\leq n}R_1^l\theta^l\Delta_{r+p}^{n-l}.$$

\noindent
To compute $O_r(n)$, in the last sum decompose $R_1^l=R_1^k(\rho_l/\rho_k)+R_{k+1}^l$. As a result~:

\begin{eqnarray}
O_r(n)&=&\sum_{1\leq k\leq n}R_1^k\sum_{-r\leq p\leq r}(-1)^{p}\theta^k\Delta_{r-p}^{n-k}\sum_{k<l\leq n}\theta^l\Delta_{r+p}^{n-l}\rho_l\nonumber\\
&+&(-1)^{r+1}\sum_{1\leq k\leq n}\theta^k\Delta_{2r+1}^{n-k}\rho_k+\sum_{1\leq k\leq n}\sum_{-r\leq p\leq r}(-1)^{p}\theta^k\Delta_{r-p}^{n-k}\theta^k\Delta_{r+p+1}^{n-k}\rho_k\nonumber\\
&=&\sum_{1\leq k\leq n}R_1^k\theta^kZ_r(n-k)\rho_k+\sum_{1\leq k\leq n}\sum_{-r\leq p\leq r-1}(-1)^{p}\theta^k\Delta_{r-p}^{n-k}\theta^k\Delta_{r+p+1}^{n-k}\rho_k.\nonumber\end{eqnarray}

\noindent
One more time, the last term is $0$. Together with \eqref{nrn} and \eqref{zmfinal} we obtain~:

 $$N_r(n)=\sum_{1\leq k\leq n}\rho_kR_1^k(\theta^kK_r(n-k)+2\theta^kZ_r(n-k))=\sum_{1\leq k\leq l\leq n}R_1^k\theta^lK_r(n-l)\rho_l2^{(k,l)}.$$
 
\noindent
This gives the announced formula and concludes the proof of the proposition.

\fin

\section{Proof of the theorem}

\subsection{Dominated variation}
For $u\in S_+^{d-1}$, the inverse functions of $n\rightarrow\varphi_{u,+}(n)$ and $n\rightarrow\varphi_u(n)$ check a dominated variation property at infinity (Feller, 1969). Notice that the latter property holds for $\psi_+^{-1}$ and $\psi_-^{-1}$, as a consequence of the first point of lemma \ref{error}.

\begin{lemme}
\label{young}

$ $

\noindent
1. For any $x\geq1$ and $K\geq 1$~:

$$\psi^{-1}(Kx)\leq 2K^2\psi^{-1}(x).$$

\noindent
2. There exists a constant $C(\delta)>0$, so that for any $u\in S_+^{d-1}$, any $x\geq1$ and $K\geq1$~:

$$\varphi_{u,+}^{-1}(Kx)\leq \frac{2K^2}{\delta}\varphi_{u,+}^{-1}(x)\mbox{ and }\varphi_u^{-1}(Kx)\leq \frac{K^2}{C(\delta)}\varphi_u^{-1}(x).$$

\end{lemme}

\noindent
\textit{Proof of the lemma~:} 

\noindent
1. Recall that $\psi^2(n)=n\({w_+(n)\circ v_+^{-1}(n)+w_-(n)\circ v_-^{-1}(n)}\)$. For $x\geq1$, let $n=\psi^{-1}(x)$, ie $\psi(n)\leq x<\psi(n+1)$. This implies that~:

$$\psi(K^2(n+1))\geq K\psi(n+1)>Kx.$$

\smallskip
\noindent
Hence $\psi^{-1}(Kx)\leq K^2(n+1)\leq2K^2n=2K^2\psi^{-1}(x)$.

\noindent
2. Let $\kappa_{u,+}(n)=\sum_{1\leq k\leq l\leq n}T_k^l(u)=\sum_{0\leq k<l\leq n}\rho_{k}\rho_l(\zeta_{k+1}^l(u))^2$, setting $\zeta_k^l(u)=\sum_{s=k}^l\eta_s.u/\rho_s$, with $\zeta_k^l(u)=0$ if $k>l$. We first claim that~:

$$\frac{\kappa_{u,+}(n)}{v_+(n)}=\frac{\kappa_{u,+}(n-1)}{v_+(n-1)}+\frac{\rho_n}{v_+(n)v_+(n-1)}\({\sum_{0\leq k< n}\rho_k\zeta_{k+1}^n(u)}\)^2.$$

\smallskip
\noindent
In particular, $n\longmapsto\kappa_{u,+}(n)/v_+(n)$ is non-decreasing. Indeed~:

$$\kappa_{u,+}(n)=\sum_{0\leq k<l\leq n}\rho_{k}\rho_l(\zeta_{k+1}^n(u))^2+\sum_{0\leq k<l\leq n}\rho_{k}\rho_l(\zeta_{l+1}^n(u))^2-2\sum_{0\leq k<l\leq n}\rho_{k}\rho_l\zeta_{k+1}^n(u)\zeta_{l+1}^n(u).$$  

\noindent
This is rewritten as~:

\begin{eqnarray}
\kappa_{u,+}(n)&=&\sum_{0\leq k<n}\rho_k(\zeta_{k+1}^n(u))^2\sum_{k<l\leq n}\rho_l+\sum_{1\leq l\leq n}\rho_l(\zeta_{l+1}^n(u))^2\sum_{0\leq k<l}\rho_k-\({\sum_{0\leq k\leq n}\rho_k\zeta_{k+1}^n(u)}\)^2\nonumber\\
&+&\sum_{0\leq k\leq n}(\rho_k)^2(\zeta_{k+1}^n(u))^2.\nonumber
\end{eqnarray}

\noindent
In other words~:

$$\kappa_{u,+}(n)=v_+(n)\sum_{0\leq k< n}\rho_k(\zeta_{k+1}^n(u))^2-\({\sum_{0\leq k\leq n}\rho_k\zeta_{k+1}^n(u)}\)^2.$$

\noindent
Next, directly from the definition of $\kappa_{u,+}(n)$, and then using the previous equality~:

$$\kappa_{u,+}(n)-\kappa_{u,+}(n-1)=\rho_n\sum_{0\leq k<n}\rho_k(\zeta_{k+1}^n(u))^2=\rho_n\frac{\kappa_{u,+}(n)+\({\sum_{0\leq k<n}\rho_k\zeta_{k+1}^n(u)}\)^2}{v_+(n)}.$$

\noindent
Observe that this is equivalent to the desired claim. 

\medskip
We next use that for all $n\geq0$, $\delta\leq \rho_{n+1}/ \rho_{n}\leq1/\delta$, hence $v_{+}(n+1)\leq (2/\delta) v_{+}(n)$. As a result $v_+\circ v_+^{-1}(n)\leq n\leq (2/\delta)v_+\circ v_+^{-1}(n)$. Hence for $x\geq 1$ and $K\geq1$~:

$$\kappa_{u,+}\circ v_+^{-1}(Kx)\geq\kappa_{u,+}\circ v_+^{-1}(x)\frac{v_+\circ v_+^{-1}(Kx)}{v_+\circ v_+^{-1}(x)}\geq \frac{\delta K}{2}\kappa_{u,+}\circ v_+^{-1}(x).$$

\smallskip
\noindent
A similar property is verified for some symmetrically defined function $\kappa_{u,-}\circ v_-^{-1}$. Notice that~:
 
$$\varphi_{u,+}^2(n)=\psi^2(n)+\kappa_{u,+}\circ v_+^{-1}(n)+\kappa_{u,-}\circ v_-^{-1}(n).$$

\medskip
\noindent
Notice that $\varphi_{u,+}(n)\rightarrow+\infty$, as $n\rightarrow+\infty$. As we showed in point one that $\psi^2(Kx)\geq K\psi^2(x)$, we obtain that for $x\geq 1$ and $K\geq1$~:

$$\varphi_{u,+}(Kx)\geq \sqrt{(\delta K/2)}\varphi_{u,+}(x).$$

\smallskip
\noindent
We conclude as in point one. Let $x\geq 1$ and $n=\varphi_{u,+}^{-1}(x)$ and $K\geq1$. Then $\varphi_{u,+}(n)\leq x<\varphi_{u,+}(n+1)$, so~:

$$\varphi_{u,+}((2K^2/\delta)(n+1))\geq K\varphi_{u,+}(n+1)>Kx.$$

\smallskip
\noindent
Consequently $\varphi_{u,+}^{-1}(Kx)\leq((2K^2)/\delta) \varphi_{u,+}^{-1}(x)$.

\medskip
It remains to show the same result for $\varphi_u$. This way, let $\kappa_u(-m,n)=\sum_{-m\leq k\leq l\leq n}T_k^l(u)$, for $m\geq1$, $n\geq1$. Then, the computation on $\kappa_{u,+}$ shows that~:

$$n\longmapsto \frac{\kappa_u(-m,n)}{(v_-(m)/a_0)+v_+(n)}\mbox{ and }m\longmapsto \frac{\kappa_u(-m,n)}{(v_-(m)/a_0)+v_+(n)}$$

\medskip
\noindent
are non-decreasing. This furnishes that for some constant $C(\delta)>0$~:

$$\kappa_{u}(-v_-^{-1}(Kx),v_+^{-1}(Kx))\geq \frac{C(\delta) K}{2}\kappa_{u}(-v_-^{-1}(x),v_+^{-1}(x)).$$

\smallskip
\noindent
As $\varphi_{u}^2(n)=\psi^2(n)+\kappa_{u}(-v_-^{-1}(n),v_+^{-1}(n))$, we conclude as before. This ends the proof of the lemma.

\fin

\subsection{Order of the real part of $1-\chi_D(ut)$}

With $u\in S_+^{d-1}$ and small $t>0$, recall the decomposition $\chi_D(ut)=\varphi_0(ut)(p'_0\chi_D^+(ut)+q'_0\chi_D^-(ut))$ and also that~:

$$\chi_D^+(ut)=f^+(ut)+O(R^+(t))\mbox{ and }\chi_D^-(ut)=f^-(ut)+O(R^-(t)),$$

\smallskip
\noindent
where the $O(~)$ are uniform in $u\in S_+^{d-1}$ and where $R^+(t)$ and $R^-(t)$ have respective orders $1/\psi_+^{-1}(1/t)$ and $1/\psi_-^{-1}(1/t)$, by lemma \ref{error}.

\begin{lemme}
\label{ouloulou}

$ $ 

\noindent
Let $R(t)=R^+(t)+R^-(t)$.

\smallskip
\noindent
1. We have $\chi_D(ut)=\varphi_0(ut)(p'_0f^+(ut)+q'_0f^-(ut))+O(R(t))$.

\smallskip
\noindent 
2. We have $t^2=O(R^+(t))$ and $t^2=O(R^-(t))$.

\smallskip
\noindent
3. We have $t\mbox{Im}(1-f^+(ut))=O(R^+(t))$ and $t\mbox{Im}(1-f^-(ut))=O(R^-(t))$.

\smallskip
\noindent
4. We have $\chi_D(ut)=(1+itm_0.ur_0/(1-r_0))(p'_0f^+(ut)+q'_0f^-(ut))+O(R(t))$ and

\begin{eqnarray}
\label{decompo}
\mbox{Re}(1-\chi_D)(ut)=p'_0\mbox{Re}(1-f^+(ut))+q'_0\mbox{Re}(1-f^-(ut))+O(R(t)).
\end{eqnarray}

\smallskip
\noindent
5. There is a constant $c>0$ so that for small $t>0$, uniformly in $u\in S_+^{d-1}$~:

$$\mbox{Re}(1-\chi_D(ut))\geq cR(t).$$

\end{lemme}

\smallskip
\noindent
\textit{Proof of the lemma~:}

\noindent
1. This follows from $\chi_D(ut)=\varphi_0(ut)(p'_0\chi_D^+(ut)+q'_0\chi_D^-(ut))$ and $\chi_D^{\pm}(ut)=f^{\pm}(ut)+O(R^{\pm}(t))$.

\medskip
\noindent
2. As $\psi_+^2(n)=nw_+\circ v_+^{-1}(n)$, for some constant $c>0$, $\psi_+(n^2)\geq cn$, so $\psi_+^{-1}(1/t)\leq c'/t^2$, $t>0$, for some constant $c'>0$. By lemma \ref{error}, $t^2=O(\psi_+^{-1}(1/t))=O(R_+(t))$, which gives the first property. The other one is proved in the same way.

\medskip
\noindent
3. We make use of proposition \ref{estim} and lemma \ref{cf}. Taking any integer $n\geq 1$ and since $f^+(ut)=A_n(ut)/B_n(ut)+O(1/v_+(n))$ (where $O(~)$ is independent on $u$ and $t$), we have~:

\begin{eqnarray}
\mbox{Im}(f^+(ut))&=&\frac{\mbox{Im}(A_n(ut))\bar{B}_n(ut))}{|B_n(ut)|^2}+O(1/v_+(n))=\frac{\sum_{0\leq r\leq n-1}t^{2r+1}N_r(n)}{\sum_{0\leq r\leq n}t^{2r}L_r(n)}+O(1/v_+(n)).\nonumber
\end{eqnarray}

\noindent
Now, see \eqref{eqlrn}, $L_r(n)=\sum_{0\leq l\leq k\leq n}\rho_k\theta^kK_r(n-k)2^{(l,k)}\rho_l\geq\sum_{1\leq k\leq n}\rho_kv_+(k)\theta^kK_r(n-k)$, where the dependence in $u\in S_+^{d-1}$ is implicit, and~:

$$N_r(n)=\sum_{1\leq k\leq l\leq n}R_1^k\theta^lK_r(n-l)\rho_l2^{(k,l)}=\sum_{1\leq k\leq n}\[{\sum_{1\leq s\leq l\leq k}\frac{\eta_s.u}{\rho_s}\rho_l2^{(l,k)}}\]\rho_k\theta^kK_r(n-k).$$

\smallskip
\noindent
As the $\eta_n.u$ are uniformly bounded by some $C/2$ (as $n$ and $u\in S_+^{d-1}$ vary), we get~:

$$|N_r(n)|\leq Cw_+(n)\sum_{1\leq k\leq n}v_+(k)\rho_k\theta^kK_r(n-k)\leq Cw_+(n)L_r(n).$$

\noindent
We finally obtain $|\mbox{Im}(f^+(ut))|\leq tw_+(n)+O(1/v_+(n))$. Let $n'=\psi_+^{-1}(1/t)$ and $n=v_+^{-1}(n')$. By definition of $\psi_+$, we have $n'w_+(n)\leq 1/t^2$. We obtain $|\mbox{Im}(f^+(ut))|\leq 1/(tn')+O(1/n')$, which is the desired result. The situation for $t|\mbox{Im}(f^-(ut))|$ is similar.

\medskip
\noindent
4. Write $\varphi_0(ut)=1+itm_0.ur_0/(1-r_0)+O(t^2)$, with $O(~)$ uniform in $u\in S_+^{d-1}$. Using the first point of the lemma, we get~:

$$\chi_D(ut)=\({1+\frac{itm_0.ur_0}{1-r_0}}\)(p'_0f^+(ut)+q'_0f^-(ut))+O(R(t)),$$

\noindent
with again an error term uniform in $u\in S_+^{d-1}$. Therefore~:

$$1-\chi_D(ut)=p'_0(1-f^+(ut))+q'_0(1-f^-(ut))-\frac{itm_0.ur_0}{1-r_0}(p'_0f^+(ut)+q'_0f^-(ut))+O(R(t)).$$

\noindent
Taking the real part~:

$$\mbox{Re}(1-\chi_D(ut))=p'_0\mbox{Re}(1-f^+(ut))+q'_0\mbox{Re}(1-f^-(ut))\nonumber\\
+\frac{tm_0.ur_0}{1-r_0}\({p'_0\mbox{Im}(f^+(ut))+q'_0\mbox{Im}(f^-(ut))}\)+O(R(t)).$$

\noindent
The third point of the lemma then gives \eqref{decompo}.

\medskip
\noindent
5. By lemma \ref{error}, for a constant $c_1>0$ independent on $u\in S_+^{d-1}$, we have for small $t>0$, $1-|\chi_D^{+}(ut)|\geq c_1R^{+}(ut)$. Idem, for some $c_2>0$, we get $1-|\chi_D^{-}(ut)|\geq c_2R^{-}(ut)$. As $\chi_D(ut)=\varphi_0(ut)(p'_0\chi_D^+(ut)+q'_0\chi_D^-(ut))$ and $|\varphi_0(ut)|\leq 1$~:

\begin{eqnarray}
\label{minoro}\mbox{Re}(1-\chi_D(ut))\geq1-|\chi_D(ut)|&\geq& 1-|p'_0\chi_D^+(ut)+q'_0\chi_D^-(ut)|\nonumber\\
&\geq& p'_0(1-|\chi_D^+(ut)|)+q'_0(1-|\chi_D^-(ut)|)\nonumber\\
&\geq &c_1p'_0R^+(t)+c_2q'_0R^-(t)\geq cR(t),\nonumber
\end{eqnarray}

\noindent
for some constant $c>0$. This completes the proof of the lemma.

\fin

\smallskip
\noindent
\begin{remark}Notice that in \cite{jb1} one always had $t=O(R^+(t))$. This is not true anymore here. For example if $\sum_{k\geq1}(1/\rho_k)<\infty$, one may check that $R^+(t)$ can have order $t^2$, as $t\rightarrow0$.
\end{remark}

\begin{prop}
\label{re}

$ $ 

\noindent
There is a constant $C\geq1$ so that for $t>0$ small enough, uniformly in $u\in S_+^{d-1}$~:

$$\frac{1}{C}\leq\varphi_{u,+}^{-1}(1/t)\mbox{Re}(1-\chi_D(ut))\leq C.$$

\end{prop}

\smallskip
\noindent
\textit{Proof of the proposition~:}

\noindent
We still fix $u\in S_+^{d-1}$ and $t>0$. Recall that $f^+(ut)=\lim_{n\rightarrow+\infty} A_n(ut)/B_n(ut)$, where $(A_n(ut))$ and $(B_n(ut))$ check proposition \ref{estim}. Fixing some $n\geq 1$, we use proposition \ref{estim} and lemma \ref{cf}~:

\begin{eqnarray}
\mbox{Re}(1-f^+(ut))&=&\mbox{Re}(1-A_n(ut)/B_n(ut))-\mbox{Re}\({\sum_{k>n}\frac{\rho_k}{B_k(ut)B_{k-1}(ut)}}\)\nonumber\\
&\leq&\frac{\mbox{Re}((B_n(ut)-A_n(ut))\bar{B}_n(ut))}{|B_n(ut)|^2}+\frac{1}{v_+(n)}\nonumber\\
&\leq&\frac{v_+(n)+\sum_{1\leq r\leq n}t^{2r}M_r(n)}{v_+(n)^2+\sum_{r=1}^nt^{2r}L_r(n)}+\frac{1}{v_+(n)}\leq\frac{1}{v_+(n)}\({2+\sum_{1\leq r\leq n}t^{2r}\frac{M_r(n)}{v_+(n)}}\),\nonumber\end{eqnarray}

\noindent
where $L_r(n)$ and $M_r(n)$ depend on $u$. By the formula for $M_r(n)$ and $K_r(n)$ in proposition \ref{estim}, $M_r(n)\leq (\sum_{1\leq k\leq l\leq n}T_k^l(u))^r2^rv_+(n)$, for $r\geq1$. Hence~:

$$M_r(v_+^{-1}(n))\leq n2^r\varphi_{u,+}^{2r}(n).$$

\smallskip
\noindent
As a result, for some constant $C>0$ independent on $u$ and any $n\geq1$~:

$$\mbox{Re}(1-f^+(ut))\leq\frac{C}{n}\[{1+\sum_{1\leq r\leq v_+^{-1}(n)}(2t^{2})^r\varphi_{u,+}^{2r}(n)}\].$$

\noindent
Choose $n=n_u(t)=\varphi_{u,+}^{-1}(1/(2t))$. In particular $\varphi_{u,+}^2(n)\leq1/(4t^2)$. We arrive at~:

$$\mbox{Re}(1-f^+(ut))\leq \frac{C}{n}\({1+\sum_{r\geq1}2^{-r}}\)\leq \frac{2C}{n}=\frac{2C}{\varphi_{u,+}^{-1}(1/(2t))}\leq \frac{C'}{\varphi_{u,+}^{-1}(1/t)},$$

\noindent
for some constant $C'$ independent on $u$, using lemma \ref{young}. Idem, $\mbox{Re}(1-f^-(ut))\leq C'/\varphi_{u,+}^{-1}(1/t)$. Via now \eqref{decompo}, using that $R^{\pm}(t)=O(1/\psi_{\pm}^{-1}(1/t))=O(1/\varphi_{u,+}^{-1}(1/t))$, this shows the right-hand side inequality of the proposition.

\medskip
Consider next the other direction. Starting in the same way, for any $n\geq1$, via proposition \ref{estim} and lemma \ref{cf} (third point)~:

\begin{eqnarray}
\mbox{Re}(1-f^+(ut))&=&\mbox{Re}(1-A_n(ut)/B_n(ut))-\mbox{Re}\({\sum_{k>n}\frac{\rho_k}{B_k(ut)B_{k-1}(ut)}}\)\nonumber\\
&\geq&\frac{\mbox{Re}((B_n(ut)-A_n(ut))\bar{B}_n(ut))}{|B_n(ut)|^2}-\frac{v_+(n)}{|B_n(ut)|^2}\nonumber\\
&=&\frac{v_+(n)+\sum_{1\leq r\leq n}t^{2r}M_r(n)}{|B_n(ut)|^2}-\frac{v_+(n)}{|B_n(ut)|^2}=\frac{\sum_{1\leq r\leq n}t^{2r}M_r(n)}{v_+(n)^2+\sum_{1\leq r\leq n}t^{2r}L_r(n)}.\nonumber
\end{eqnarray}

\noindent
By prop. \ref{estim}, $M_r(n)=\sum_{0\leq k\leq n}\rho_k\theta^kK_r(n-k)$ and $L_r(n)=\sum_{0\leq l\leq k\leq n}2^{(l,k)}\rho_l\rho_k\theta^kK_r(n-k)$, so we have $L_r(n)\leq 2v_+(n)M_r(n)$. Hence~:

\begin{equation}
\label{yarrive}\mbox{Re}(1-f^+(ut))\geq \frac{1}{v_+(n)}\frac{\sum_{1\leq r\leq n}t^{2r}M_r(n)/v_+(n)}{1+2\sum_{1\leq r\leq n}M_r(n)/v_+(n)}\geq \frac{1}{v_+(n)}\frac{t^2M_1(n)/v_+(n)}{1+2t^2M_1(n)/v_+(n)},
\end{equation}

\smallskip
\noindent
using in the last step that $x\longmapsto x/(1+2x)$ is increasing ($x>0$). As a result, for some constant $c>0$ independent on $u$ and all $n\geq1$~:

$$\mbox{Re}(1-f^+(ut))\geq\frac{c}{n}\frac{ct^2M_1(v_+^{-1}(n))/n}{1+2ct^2M_1(v_+^{-1}(n))/n}.$$

\medskip
Let $\kappa_{u,+}(m)=\sum_{1\leq k\leq l\leq m}T_k^l(u)$ and assume first that $\lim_{m\rightarrow+\infty}\kappa_{u,+}(m)=+\infty$. Note (using proposition \ref{estim}) that $M_1(n)\geq \sum_{1\leq m\leq n}\rho_m\kappa_{u,+}(m)$. Therefore~:

$$M_1(v_+^{-1}(n))\geq \sum_{1\leq m\leq v_+^{-1}(n)}\rho_m\kappa_{u,+}(m).$$

\noindent
Let $c_0\geq2$ be such that for all $n$, $v_+(n+1)\leq c_0v_+(n)$. Set $m_u(t)=(\kappa_{u,+}\circ v_+^{-1})^{-1}(1/t^2)$ and next choose $n_u(t)=c_0^2m_u(t)$. Let $s=v_+^{-1}(m_u(t))$ and $s'=v_+^{-1}(n_u(t))$. This gives~:

$$v_+(s)\leq m_u(t)<v_+(s+1)\leq c_0v_+(s)\mbox{ and }v_+(s')\leq c_0^2m_u(t)<v_+(s'+1)\leq c_0v_+(s').$$

\smallskip
\noindent
As a result, $c_0^2m_u(t)\geq v_+(s')-v_+(s)\geq (c_0-1)m_u(t)\mbox{ and }m_u(t)\geq v_+(s)\geq m_u(t)/c_0$. This furnishes the inequalities~:

$$\frac{M_1(v_+^{-1}(n_u(t)))}{n_u(t)}\geq\frac{\sum_{s<m\leq s'}\rho_m\kappa_+(m)}{n_u(t)}\geq \kappa_{u,+}(s+1)\frac{v_+(s')-v_+(s)}{n_u(t)}\geq\frac{\alpha}{t^2}\mbox{ with } \alpha=(c_0-1)/c_0^2.$$

\smallskip
\noindent
Consequently, with $\alpha'=(c^2\alpha)/(c_0^2(1+2c\alpha))$~:

$$\mbox{Re}(1-f^+(ut))\geq\frac{c}{n_u(t)}\frac{c\alpha}{1+2c\alpha}=\frac{\alpha'}{(\kappa_{u,+}\circ v_+^{-1})^{-1}(1/t^2)}.$$

\smallskip
If now $m\longmapsto \kappa_{u,+}(m)$ is bounded, the previous inequality is valid as long as $(\kappa_{u,+}\circ v_+^{-1})^{-1}(1/t^2)$ is defined. For smaller $t$, we have $(\kappa_{u,+}\circ v_+^{-1})^{-1}(1/t^2)=+\infty$, so that the previous lower-bound is obvious in this case. Similarly, with $\kappa_{u,-}(m)=\sum_{-m\leq k\leq l\leq -1}T_k^l(u)$, we have~:

$$\mbox{Re}(1-f^-(ut))\geq\alpha'/(\kappa_{u,-}\circ v_-^{-1})^{-1}(1/t^2).$$

\smallskip
\noindent
To prove a lower bound, we use \eqref{decompo}, giving for some constant $c_3>0$ independent on $u$~:

\begin{equation}
\label{minno}
\mbox{Re}(1-\chi_D(ut))\geq p'_0\mbox{Re}(1-f^+(ut))+q'_0\mbox{Re}(1-f^-(ut))-c_3/\psi^{-1}(1/t).
\end{equation}

\medskip
\noindent
Recall that $\varphi_{u,+}^2=\psi^2+\kappa_{u,+}\circ v_+^{-1}+\kappa_{u,-}\circ v_-^{-1}$. Then, for some constant $\beta>0$ independent on $u$ and $t$, we have~:

$$\beta\leq\frac{\varphi_{u,+}^{-1}(1/t)}{\min\{\psi^{-1}(1/t), (\kappa_{u,+}\circ v_+^{-1})^{-1}(1/t^2),(\kappa_{u,-}\circ v_-^{-1})^{-1}(1/t^2)\}}\leq1.$$

\smallskip
\noindent
Fixing $t>0$, suppose for example that $(\kappa_{u,+}\circ v_+^{-1})^{-1}(1/t^2)\leq (\kappa_{u,-}\circ v_-^{-1})^{-1}(1/t^2)$. This leads to the following discussion~:

\medskip

\noindent
-- If $(1/\psi^{-1}(1/t))\leq p_0'\mbox{Re}(1-f^+(ut))/(2c_3)$ and $\psi^{-1}(1/t)\geq(\kappa_{u,+}\circ v_+^{-1})^{-1}(1/t^2)$, then~:

$$\mbox{Re}(1-\chi_D(ut))\geq (p'_0/2)\mbox{Re}(1-f^+(ut))\geq\frac{p'_0\alpha'/2}{(\kappa_{u,+}\circ v_+^{-1})^{-1}(1/t^2)}\geq\frac{\beta p'_0\alpha'/2}{\varphi_{u,+}^{-1}(1/t)}.$$

\medskip
\noindent
-- If $(1/\psi^{-1}(1/t))\leq p_0'\mbox{Re}(1-f^+(ut))/(2c_3)$ and $\psi^{-1}(1/t)\leq(\kappa_{u,+}\circ v_+^{-1})^{-1}(1/t^2)$, then, by lemma \ref{ouloulou} and proposition \ref{error}, for absolute constants $c>0$ and $c'>0$~:

$$\mbox{Re}(1-\chi_D(ut))\geq cR(ut)\geq c'/\psi^{-1}(1/t)\geq \beta c'/\varphi_{u,+}^{-1}(1/t).$$

\medskip
\noindent
-- If $(1/\psi^{-1}(1/t))>p_0'\mbox{Re}(1-f^+(ut))/(2c_3)$, then $\psi^{-1}(1/t)<(2c_3/(p'_0\alpha'))(\kappa_{u,+}\circ v_+^{-1})^{-1}(1/t^2)$. We obtain the inequality~:

$$\frac{\beta}{\varphi_{u,+}^{-1}(1/t)}\leq \frac{1}{\psi^{-1}(1/t)\min\{p'_0\alpha'/(2c_3),1\}}.$$

\smallskip
\noindent
We conclude as in the previous case, via $\mbox{Re}(1-\chi_D(ut))\geq cR(ut)\geq c'/\psi^{-1}(1/t)$. This completes the proof of the proposition.

\fin

\subsection{Preliminaries for estimating the modulus of $1-\chi_D(ut)$}

We still fix $u\in S_+^{d-1}$ and $t>0$. We use proposition \ref{estim} concerning $f^+$ and its symmetric analogue for $f^-$. To precise the dependency with respect to $f^+$ or $f^-$, we put a superscript ($+$ or $-$) on $A_n$, $B_n$, etc. For example $f^{+}(ut)=\lim_{n\rightarrow+\infty}A_n^{+}(ut)/B_n^{+}(ut)$. Keeping the same sets of summation, the expressions corresponding to $K_r^-(n)$, etc, are deduced from proposition \ref{estim} by replacing $(q_k,p_k)$ by $(p_{-k},q_{-k})$. Any $\rho_k$ becomes $\rho_{-k-1}q_0/p_0$. It is worth noticing that $T_k^l(u)$ is simply transformed into $T_{-l}^{-k}(u)$.

\medskip
Let us begin with a formal computation on reversed continued fractions.

\begin{lemme}
\label{reverse}

$ $

\noindent
Let $n\geq1$ and consider the formal reduced continued fraction~:

$$\frac{U_n}{V_n}=[(-c_1,d_1);(-c_2,d_2);\cdots;(-c_n,d_n)].$$

\noindent
Then the reduced reversed continued fraction~:

$$\frac{\tilde{U_n}}{\tilde{V_n}}=[(-1/c_n,d_n/c_n);(-1/c_{n-1},d_{n-1}/c_{n-1});\cdots;(-1/c_1,d_1/c_1)]$$

\noindent
verifies $V_n=c_1\cdots c_n\tilde{V_n}$.

\end{lemme}

\smallskip
\noindent
\textit{Proof of the lemma~:}

\noindent
We have~:

$$V_n=\langle e_1,\begin{pmatrix}
d_n&-c_n\\
1&0\\
\end{pmatrix}\cdots\begin{pmatrix}
d_1&-c_1\\
1&0\\
\end{pmatrix} 
e_1\rangle$$

\noindent
Transposing and next conjugating the matrices with $\mbox{diag}(1,-1)$~:

\begin{eqnarray}
V_n=\langle e_1,\begin{pmatrix}
d_1&1\\
-c_1&0\\
\end{pmatrix}\cdots\begin{pmatrix}
d_n&1\\
-c_n&0\\
\end{pmatrix} 
e_1\rangle&=&\langle e_1,\begin{pmatrix}
d_1&-1\\
c_1&0\\
\end{pmatrix}\cdots\begin{pmatrix}
d_n&-1\\
c_n&0\\
\end{pmatrix} 
e_1\rangle\nonumber\\
&=&c_1\cdots c_n\langle e_1,\begin{pmatrix}
d_1/c_1&-1/c_1\\
1&0\\
\end{pmatrix}\cdots\begin{pmatrix}
d_n/c_n&-1/c_n\\
1&0\\
\end{pmatrix} 
e_1\rangle.\nonumber\nonumber\end{eqnarray}

\noindent
Hence $V_n=c_1\cdots c_n \tilde{V_n}$. This proves the lemma.

\fin

\medskip

Let us start from relation \eqref{gw}, $\chi_D(ut)=\varphi_0(ut)(p'_0\chi_D^+(ut)+q'_0\chi_D^-(ut))=(\varphi_0(ut)/b_0)(\chi_D^+(ut)+a_0\chi_D^-(ut))$. This gives, using lemmas \ref{error} and \ref{ouloulou} and taking $t>0$ small, independently on $u$~:

\begin{eqnarray}\chi_D(ut)-1&=&\frac{\varphi_0(ut)}{b_0}(\chi_D^+(ut)+a_0\chi_D^-(ut)-b_0/\varphi_0(ut))\nonumber\\
&=&\frac{\varphi_0(ut)}{b_0}(f^+(ut)+a_0f^-(ut)-b_0/\psi_0(ut))+O(R(t))\nonumber\\
&=&\frac{\varphi_0(ut)}{b_0}\({\frac{A_n^+(ut)}{B_n^+(ut)}+a_0\frac{A_m^-(ut)}{B_m^-(ut)}-b_0/\psi_0(ut)}\)\nonumber\\
&+&\frac{\varphi_0(ut)}{b_0}\({\sum_{k>n}\frac{\rho_k}{B_k^+(ut)B_{k-1}^+(ut)}+\sum_{k>m}\frac{a_0^2\rho_{-k-1}}{B_k^-(ut)B_{k-1}^-(ut)}}\)+O(R(t)),\nonumber\end{eqnarray}

\noindent
with $O(~)$ uniform in $u$ and arbitrary $n\geq1$, $m\geq1$. As a result~:

\begin{eqnarray}\chi_D(ut)-1&=&\frac{\varphi_0(ut)}{b_0B_n^+(ut)B_m^-(ut)}\({A_n^+(ut)B_m^-(ut)+a_0A_m^-(ut)B_n^+(ut)-(b_0/\psi_0(ut))B_n^+(ut)B_m^-(ut)}\)\nonumber\\
&+&\frac{\varphi_0(ut)}{b_0}R_{-m,n}(ut)+O(R(t)),
\end{eqnarray}

\noindent
with $|R_{-m,n}(ut)|\leq (v_+(n)/|B_n^+(ut)|^2)+a_0(v_-(m)/|B_m^-(ut)|^2)$, by proposition \ref{cf}, and $O(~)$ uniform in $u\in S_+^{d-1}$.

\begin{lemme}

$ $

\noindent
Let $n\geq1$, $m\geq1$ and the following reduced continued fraction~:

$$\frac{\tilde{A}_{m+n+1}(ut)}{\tilde{B}_{m+n+1}(ut)}=[(-a_{-m},b_{-m}/\psi_{-m}(ut));(-a_{-m+1},b_{-m+1}/\psi_{-m+1}(ut));\cdots;(-a_{n},b_{n}/\psi_{n}(ut))].$$

\smallskip
\noindent
Then $\tilde{B}_{m+n+1}(ut)=-a_{-1}\cdots a_{-m}(A_n^+(ut)B_m^-(ut)+a_0A_m^-(ut)B_n^+(ut)-(b_0/\psi_0(ut))B_n^+(ut)B_m^-(ut))$.
\end{lemme}

\medskip
\noindent
\textit{Proof of the lemma~:}

\noindent
Fix $m\geq1$. Observe that the two functions $n\longmapsto -\tilde{B}_{m+n+1}(ut)/(a_{-1}\cdots a_{-m})$ and $n\longmapsto A_n^+(ut)B_m^-(ut)+a_0A_m^-(ut)B_n^+(ut)-(b_0/\psi_0(ut))B_n^+(ut)B_m^-(ut)$ check the same recursive relation $X_n=(b_n/\psi_n(ut))X_{n-1}-a_nX_{n-2}$, for $n\geq1$. We just need to check that they coincide for the values $n=0$ and $n=1$.

\medskip
First, $\tilde{B}_m(ut)/(a_{-1}\cdots a_{-m})=B_m^-(ut)$ and~:

$$A_m^-(ut)=(1/a_{-1})\theta^{-1}B^-_{m-1}(ut)=\tilde{B}_{m-1}(ut)/(a_{-1}\cdots a_{-m}),$$

\smallskip
\noindent
by lemma \ref{reverse}. For $n=0$ we have $-\tilde{B}_{m+1}(ut)/(a_{-1}\cdots a_{-m})$ and $a_0A_m^-(ut)-(b_0/\psi_0(ut))B_m^-(ut)$. Since one has~:

$$\tilde{B}_{m+1}(ut)=(b_0/\psi_0(ut))\tilde{B}_m(ut)-a_0\tilde{B}_{m-1}(ut),$$

\smallskip
\noindent
this gives the result for $n=0$. For $n=1$, we have~:

$$\tilde{B}_{m+2}(ut)=\frac{b_1}{\psi_1(ut)}\tilde{B}_{m+1}(ut)-a_1\tilde{B}_m(ut)=\({\frac{b_1}{\psi_1(ut)}\frac{b_0}{\psi_0(ut)}-a_1}\)\tilde{B}_m(ut)-\frac{b_1}{\psi_1(ut)}a_0\tilde{B}_{m-1}(ut).$$

\smallskip
\noindent
This has to be compared with $a_1B_m^-(ut)+a_0(b_1/\psi_1(ut))A_m^-(ut)-(b_0/\psi_0(ut))(b_1/\psi_1(ut))B_m^-(ut)$. This provides the conclusion of the lemma.\fin

\bigskip
As a consequence of this lemma we obtain~:

\begin{equation}
\label{bondebut}
\chi_D(ut)-1=-\frac{\varphi_0(ut)}{b_0}\({a_0\rho_{-m-1}\frac{\tilde{B}_{m+n+1}(ut)}{B_n^+(ut)B_m^-(ut)}-R_{-m,n}(ut)}\)+O(R(t)).
\end{equation}

\medskip
\noindent
Now it follows from proposition \ref{estim} that $|B^+_n(ut)|^2=\sum_{r=0}^nt^{2r}L^+_r(n)$, with $H_r((k_i),(l_j)):=\#\{0\leq i\leq r~|~l_i+1<k_{i+1}\}$ and~:

$$L^+_r(n)=\sum_{0\leq l_0<k_1\leq l_1<\cdots<k_r\leq l_r< k_{r+1}\leq n+1}\rho_{l_0}T_{k_1}^{l_1}(u)\cdots T_{k_{r}}^{l_{r}}(u)\rho_{k_{r+1}-1}2^{H_r((k_i),(l_j))}.$$

\medskip
\noindent
As a result, setting $W_{-m,n}(ut)=a_0\rho_{-m-1}\tilde{B}_{m+n+1}(ut)$, we have~:

$$|W_{-m,n}(ut)|^2=\sum_{0\leq r\leq n+m+1}t^{2r}U_r(u),$$

\medskip
\noindent
with $H_r((k_i),(l_i))=\#\{0\leq i\leq r~|~l_i+1<k_{i+1}\}$ and~:

\begin{eqnarray}
U_r(u)=a_0^2\rho_{-m-1}^2&\sum_{-m-1\leq l_0<k_1\leq l_1<\cdots<k_s\leq l_s<k_{s+1}\leq n+1}\theta^{-m-1}&\rho_{l_0+m+1}T_{k_1}^{l_1}(u)\cdots T_{k_s}^{l_s}(u)\theta^{-m-1}\nonumber\\
&&\times\rho_{k_{s+1}-1+m+1}2^{H_s((k_i),(l_i))}.\nonumber
\end{eqnarray}

\smallskip
\noindent
After a cocycle simplification~:

\begin{equation}
\label{usa}U_r(u)=a_0^2\sum_{-m-1\leq l_0<k_1\leq l_1<\cdots<k_s\leq l_s<k_{s+1}\leq n+1}\rho_{l_0}T_{k_1}^{l_1}(u)\cdots T_{k_s}^{l_s}(u)\rho_{k_{s+1}-1}2^{H_s((k_i),(l_i))}.
\end{equation}

\subsection{Order of the modulus of $1-\chi_D(ut)$}

\begin{prop}
\label{mod}

$ $ 

\noindent
There is a constant $C\geq1$ so that for $t>0$ small enough, uniformly in $u\in S_+^{d-1}$~:

$$\frac{1}{C}\leq\varphi_u^{-1}(1/t)|1-\chi_D(ut)|\leq C.$$

\end{prop}

\medskip
\noindent
\textit{Proof of the proposition~:}

\noindent
Let us start from \eqref{bondebut}. Set~:

$$F(ut)=\frac{W_{-m,n}(ut)}{B_n^+(ut)B_m^-(ut)}-R_{-m,n}(ut).$$

\medskip
\noindent
We then have $\chi_D(ut)-1=-(\varphi_0(ut)/b_0)F(ut)+O(R(t))$, where the last term is uniform in $u\in S_+^{d-1}$. By definition, $F(ut)$ does not depend on $m\geq1$, $n\geq1$ (see \eqref{bondebut}). The latter are arbitrary for the moment. Using the upper-bound on $R_{-m,n}(ut)$, we get~:

\begin{eqnarray}
|F(ut)|&\geq&\frac{|W_{-m,n}(ut)|}{|B_n^+(ut)||B_m^-(ut)|}-\frac{v_+(n)}{|B_n^+(ut)|^2}-a_0\frac{v_-(m)}{|B_m^-(ut)|^2} \nonumber\\
&\geq &\frac{1}{|B_n^+(ut)||B_m^-(ut)|}\({|W_{-m,n}(ut)|-v_+(n)\frac{|B_m^-(ut)|}{|B_n^+(ut)|}-a_0v_-(m)\frac{|B_n^+(ut)|}{|B_m^-(ut)|}}\).\nonumber
\end{eqnarray}

\medskip
\noindent
Recall that $W^2_{-m,n}(ut)=\sum_{0\leq r\leq n+m+1}t^{2r}U_r(u)$, with $U_s(u)$ given by \eqref{usa}. In particular~:

$$U_0(u)=a_0^2\({\sum_{-m-1\leq l\leq n}\rho_l}\)^2=(v_-(m)+a_0v_+(n))^2.$$

\noindent
Introduce $Z^2_{-m,n}(t)$ such that~:

$$|W_{-m,n}(ut)|^2-(a_0^2|B_n^+(ut)|^2+|B_m^-(ut)|^2)-2a_0v_+(n)v_-(m)=Z_{-m,n}^2(ut).$$

\medskip
\noindent
Then $Z_{-m,n}^2(ut)=\sum_{1\leq s\leq n+m+1}t^{2s}V_s(u)$, where~:

\begin{equation}
\label{etvs}
V_s(u)=a_0^2\sum_{\substack{-m-1\leq l_0<k_1\leq l_1<\cdots<k_s\leq l_s<k_{s+1}\leq n+1\\l_0<0<k_{s+1}}}\rho_{l_0}T_{k_1}^{l_1}(u)\cdots T_{k_s}^{l_s}(u)\rho_{k_{s+1}-1}2^{H_s((k_i),(l_i))}.
\end{equation}

\medskip
\noindent
Observe now that~:

\begin{eqnarray}
\({v_+(n)\frac{|B_m^-(ut)|}{|B_n^+(ut)|}+a_0v_-(m)\frac{|B_n^+(ut)|}{|B_m^-(ut)|}}\)^2&=&v_+^2(n)\frac{|B_m^-(ut)|^2}{|B_n^+(ut)|^2}+a^2_0v_-^2(m)\frac{|B_n^+(ut)|^2}{|B_m^-(ut)|^2}+2a_0v_+(n)v_-(m)\nonumber\\
&\leq&|B_m^-(ut)|^2+a^2_0|B_n^+(ut)|^2+2a_0v_+(n)v_-(m)\nonumber\\
&\leq&|W_{-m,n}(ut)|^2-Z_{-m,n}^2(ut)\leq |W_{-m,n}(ut)|^2.\end{eqnarray}

\noindent
This allows to write~:

\begin{eqnarray}
\label{debmino}
|F(ut)|&\geq&\frac{1}{|B_n^+(ut)||B_m^-(ut)|}\({|W_{-m,n}(ut)|-v_+(n)\frac{|B_m^-(ut)|}{|B_n^+(ut)|}-a_0v_-(m)\frac{|B_n^+(ut)|}{|B_m^-(ut)|}}\) \nonumber\\
&\geq &\frac{|W_{-m,n}(ut)|^2-\({v_+(n)\frac{|B_m^-(ut)|}{|B_n^+(ut)|}+a_0v_-(m)\frac{|B_n^+(ut)|}{|B_m^-(ut)|}}\)^2}{2W_{-m,n}(ut)|B_n^+(ut)||B_m^-(ut)|}\nonumber\\
&\geq&\frac{Z^2_{-m,n}(ut)}{2|W_{-m,n}(ut)||B_n^+(ut)||B_m^-(ut)|}.\end{eqnarray}

\medskip
\noindent
We now give upper-bounds on $|W_{-m,n}(ut)|$ and $|B_n^+(ut)||B_m^-(ut)|$. Observe first that $L^+_r(n)\leq V_r(u)v_+(n)/(a_0v_-(m))$, for $r\geq 1$, so that~:

$$|B_n^+(ut)|^2-v_+(n)^2\leq\frac{v_+(n)}{a_0v_-(m)}.$$

\noindent
Similarly, $|B_m^-(ut)|^2-v_-(m)^2\leq Z^2_{-m,n}(ut)v_-(m)/(a_0v_+(n))$. We obtain~:

\begin{eqnarray}
\label{minow}|W_{-m,n}(ut)|^2&=&(v_-(m)+a_0v_+(n))^2+a_0^2(|B^+_n(ut)|^2-v_+(n)^2)+(|B_m^-(ut)|^2-v_-(m)^2)+Z^2_{-m,n}(ut)\nonumber\\
&\leq&(v_-(m)+a_0v_+(n))^2+(a_0v_+(n)/v_-(m)+v_-(m)/(a_0v_+(n))+1)Z^2_{-m,n}(ut)\nonumber\\
&\leq&(v_-(m)+a_0v_+(n))^2\[{1+\frac{1}{a_0v_-(m)v_+(n)}Z^2_{-m,n}(ut)}\].\end{eqnarray}

\medskip
\noindent
In the same way~:

\begin{eqnarray}
|B_n^+(ut)|^2|B_m^-(ut)|^2&=&v_+(n)^2v_-(m)^2\({1+\sum_{1\leq r\leq n}t^{2r}\frac{L_r^+(n)}{v_+(n)^2}}\)\({1+\sum_{1\leq r\leq m}t^{2r}\frac{L_r^-(n)}{v_-(m)^2}}\)\nonumber\\
&=&v_+(n)^2v_-(m)^2\({1+\sum_{1\leq s\leq m+n}t^{2s}\sum_{0\leq r\leq s}\frac{L_r^+(n)L_{s-r}^-(m)}{v_+(n)^2v_-(m)^2}}\).\nonumber\end{eqnarray}

\medskip
\noindent
Notice that $\sum_{0\leq r\leq s}L_r^+(n)L_{s-r}^-(m)\leq V_s(u)v_+(n)v_-(m)/a_0$. Therefore~:

$$|B_n^+(ut)|^2|B_m^-(ut)|^2\leq v_+(n)^2v_-(m)^2\({1+\frac{Z^2_{-m,n}(ut)}{a_0v_+(n)v_-(m)}}\)$$

\medskip
\noindent
Inserting these two upper-bounds in \eqref{debmino} and using in the last step that the function $x\longmapsto/(1+x)$ is increasing, we obtain~:

\begin{eqnarray}
|F(ut)|&\geq&\frac{1}{2(v_-(m)/a_0+v_+(n))}\frac{Z^2_{-m,n}(ut)/(a_0v_+(n)v_-(m))}{1+Z^2_{-m,n}(ut)/(a_0v_+(n)v_-(m))}\nonumber\\
&\geq&\frac{1}{2(v_-(m)/a_0+v_+(n))}\frac{t^2V_1(u)/(a_0v_+(n)v_-(m))}{1+t^2V_1(u)/(a_0v_+(n)v_-(m))}.\nonumber
\end{eqnarray}

\medskip
\noindent
Let us now focus on $V_1(u)$ that we write $V_1(u)=V_{u,1}(-m,n)$. Set $\kappa_u(r,s)=\sum_{r\leq k\leq l\leq s}T_k^l(u)$, for $r\leq s$. We assume first that $\sum_{-\infty<k\leq l<+\infty}T_k^l(u)=+\infty$. We obtain~:

\begin{eqnarray}
V_{u,1}(-m,n)&=&a^2_0\sum_{\substack{-m-1\leq l_0<k_1\leq l_1<k_{2}\leq n+1\\l_0<0<k_{2}}}\rho_{l_0}T_{k_1}^{l_1}(u)\rho_{k_{2}-1}2^{H_1((k_i),(l_i))}\nonumber\\
&\geq&a_0^2\sum_{-m\leq l_0\leq0\leq k_2\leq n}\rho_{l_0-1}\rho_{k_2}\kappa_u(l_0,k_2).\nonumber
\end{eqnarray}

\medskip
\noindent
We next have the existence of a constant $c>0$ independent on $u\in S_+^{d-1}$ so that for all $n\geq1$~:

$$|F(ut)|\geq \frac{c}{n}\frac{(ct^2/n^2)V_{u,1}(-v_-^{-1}(n),v_+^{-1}(n))}{1+(ct^2/n^2)V_{u,1}(-v_-^{-1}(n),v_+^{-1}(n))}.$$

\smallskip
\noindent
Let $c_0\geq2$ be such that for all $n\geq 0$, $v_+(n+1)\leq c_0v_+(n)$ and $v_-(n+1)\leq c_0v_-(n)$. Taking $t>0$, set $m_u(t)=\kappa_u(-v_-^{-1}(.),v_+^{-1}(.))^{-1}(1/t^2)$. Choose next $n_u(t)=c_0^2m_u(t)$. Let $r=v_+^{-1}(m_u(t))$, $s=v_-^{-1}(m_u(t))$ and $r'=v_+^{-1}(n_u(t))$, $s'=v_-^{-1}(n_u(t))$. This gives~:

$$v_+(r)\leq m_u(t)<v_+(r+1)\leq c_0v_+(r)\mbox{ and }v_+(r')\leq c_0^2m_u(t)<v_+(r'+1)\leq c_0v_+(r').$$

\medskip
\noindent
As a result, $c_0^2m_u(t)\geq v_+(r')-v_+(r)\geq (c_0-1)m_u(t)\mbox{ and }m_u(t)\geq v_+(r)\geq m_u(t)/c_0$. In the same way, we have~:

$$v_-(s)\leq m_u(t)<v_-(s+1)\leq c_0v_-(s)\mbox{ and }v_-(s')\leq c_0^2m_u(t)<v_+(s'+1)\leq c_0v_+(s').$$

\smallskip
\noindent
Also, $c_0^2m_u(t)\geq v_-(s')-v_-(s)\geq (c_0-1)m_u(t)\mbox{ and }m_u(t)\geq v_-(s)\geq m_u(t)/c_0$.
\smallskip
\noindent
We obtain~:

\begin{eqnarray}
\frac{V_{u,1}(-v_-^{-1}(n_u(t)),v_+^{-1}(n_u(t)))}{n_u(t)^2}&\geq& a^2_0\frac{\sum_{r<l\leq r',s< k\leq s'}\rho_{-k-1}\rho_l\kappa_u(-k,l)}{n_u(t)^2}\nonumber\\
&\geq&a_0\kappa_u(-s-1,r+1)\frac{(v_+(r')-v_+(r))(v_-(s')-v_-(s))}{n_u(t)^2}\nonumber\\
&\geq&\frac{a_0(c_0-1)^2m_u(t)^2}{t^2n_u(t)^2}=\frac{\alpha}{t^2},\nonumber\end{eqnarray}

\noindent
where $\alpha=a_0(c_0-1)^2/c_0^2$. The conclusion for the moment is that there is a constant $c'>0$ independent on $u\in S_+^{d-1}$ so that for small $t>0$~:

$$|F(ut)|\geq \frac{c'}{\kappa_u(-v_-^{-1}(.),v_+^{-1}(.))^{-1}(1/t^2)}.$$

\noindent
When $\sum_{-\infty<k\leq l<+\infty}T_k^l$ is bounded, the inequality is verified, as $\kappa_u(v_-^{-1}(.),v_+^{-1}(.))^{-1}(1/t^2)=+\infty$, for small enough $t>0$. The previous lower-bound is then obvious in that case.

\medskip
In order to draw the conclusion, recall that $\varphi_u^2(n)=\psi^2(n)+\kappa_u(-v_-^{-1}(n),v_+^{-1}(n))$ and $1-\chi_D(ut)=(\varphi_0(ut)/b_0)F(ut)+O(R(t))$, with $O(~)$ uniform in $u\in S_+^{d-1}$. Also, by lemma \ref{ouloulou}~:

$$|1-\chi_D(ut)|\geq\mbox{Re}(1-\chi_D(ut))\geq c_1R(ut),$$

\smallskip
\noindent
for some absolute constant $c_1>0$. Similarly, for constants $c_2>0$ and $c_3>0$, we have the inequalities $c_2\leq R(t)\psi^{-1}(1/t)\leq c_3$. Then, for constants $\beta>0$ and $c_4>0$ independent on $u\in S_+^{d-1}$, for small $t>0$~:

$$\beta\leq\frac{\varphi_u^{-1}(1/t)}{\min\{\psi^{-1}(1/t),\kappa_u(-v_-^{-1}(.),v_+^{-1}(.))^{-1}(1/t^2)\}}$$

\noindent
and~:

$$|1-\chi_D(ut)|\geq\frac{1}{2b_0}|F(ut)|-c_4R(t)\geq \frac{c'}{2b_0\kappa_u(-v_-^{-1}(.),v_+^{-1}(.))^{-1}(1/t^2)}-\frac{c_4c_3}{\psi^{-1}(1/t)}.$$

\smallskip
\noindent
Fixing $t>0$, we then have the following discussion~:

\smallskip
\noindent
-- If $c'/(2b_0\kappa_u(-v_-^{-1}(.),v_+^{-1}(.))^{-1}(1/t^2))\geq 2c_3c_4/\psi^{-1}(1/t)$ and $\kappa_u(-v_-^{-1}(.),v_+^{-1}(.))^{-1}(1/t^2)\leq \psi^{-1}(1/t)$~:

$$|1-\chi_D(ut)|\geq \frac{c'}{4b_0\kappa_u(-v_-^{-1}(.),v_+^{-1}(.))^{-1}(1/t^2)}\geq \frac{c'\beta}{4b_0\varphi_u^{-1}(1/t)}.$$

\smallskip
\noindent
-- If $c'/(2b_0\kappa_u(-v_-^{-1}(.),v_+^{-1}(.))^{-1}(1/t^2))\geq 2c_3c_4/\psi^{-1}(1/t)$ and $\kappa_u(-v_-^{-1}(.),v_+^{-1}(.))^{-1}(1/t^2)>\psi^{-1}(1/t)$~:

$$|1-\chi_D(ut)|\geq c_1c_2/\psi^{-1}(1/t)\geq c_1c_2\beta/\varphi_u^{-1}(1/t).$$

\smallskip
\noindent
-- If $c'/(2b_0\kappa_u(-v_-^{-1}(.),v_+^{-1}(.))^{-1}(1/t^2))<2c_3c_4/\psi^{-1}(1/t)$, then for some absolute constant $c_5>0$ (independent on $u$), $1/\psi^{-1}(1/t)\geq c_5/\varphi_u^{-1}(1/t)$. Then, as above~:

$$|1-\chi_D(ut)|\geq c_1c_2/\psi^{-1}(1/t)\geq c_1c_2c_5\varphi_u^{-1}(1/t).$$

\medskip
This completes the proof of the lower bound. We next turn to the proof of the upper-bound. Let us start from the following inequality, for any $m\geq1$, $n\geq1$, using lemma \ref{cf}~:

\begin{eqnarray}
|1-\chi_D(ut)|&\leq&\frac{1}{b_0}|F(ut)|+O(R(t))\leq\frac{|W_{-m,n}(ut)|}{b_0|B_n^+(ut)||B_m^-(ut)|}+\frac{v_+(n)}{|B_n^+(ut)|^2}+\frac{a_0v_-(m)}{|B_m^-(ut)|^2}+O(R(t)),\nonumber
\end{eqnarray}

\noindent
with $O(~)$ uniform in $u\in S_+^{d-1}$. Observe that from the second line in \eqref{debmino}~:

$$\frac{v_+(n)}{|B_n^+(ut)|^2}+\frac{a_0v_-(m)}{|B_m^-(ut)|^2}\leq \frac{1}{|B_n^+(ut)||B_m^-(ut)|}\({v_+(n)\frac{|B_m^-(ut)|}{|B_n^+(ut)|}+a_0v_-(m)\frac{|B_n^+(ut)|}{|B_m^-(ut)|}}\)\leq \frac{|W_{-m,n}(ut)|}{|B_n^+(ut)||B_n^-(ut)|}.$$

\medskip
\noindent
Since $R(t)=O(1/\psi^{-1}(1/t))=O(1/\varphi_u^{-1}(1/t))$, uniformly on $u\in S_+^{d-1}$, there exists some absolute constant $C>0$ such that for small $t>0$ and all $m\geq1$ and $n\geq1$~:

$$|1-\chi_D(ut)|\leq C\frac{|W_{-m,n}(ut)|}{|B_n^+(ut)||B_m^-(ut)|}+\frac{C}{\varphi_u^{-1}(1/t)}.$$

\medskip
\noindent
From \eqref{minow} and lemma \ref{cf}, we have~:

$$\frac{|W_{-m,n}(ut)|}{|B_n^+(ut)||B_m^-(ut)|}\leq\frac{(v_-(m)+a_0v_+(n))\sqrt{1+Z^2_{-m,n}(ut)/(a_0v_-(m)v_+(n))}}{v_+(n)v_-(m)}.$$

\medskip
\noindent
Let us recall that $Z^2_{-m,n}(ut)=\sum_{1\leq s\leq m+n+1}t^{2s}V_s(u)$, where $V_s(u)$ is given by relation \eqref{etvs}, so checks $V_s(u)\leq a_0v_-(m)v_+(n)\kappa_u(-m,n)^s$, still setting $\kappa_u(-m,n)=\sum_{-m\leq k\leq l\leq n}T_k^l(u)$. As a result, for another constant $C>0$ independent on $u\in S_+^{d-1}$, small $t>0$ and any $n\geq1$~:

\begin{eqnarray}
|1-\chi_D(ut)|&\leq& \frac{C}{n}\sqrt{1+\sum_{1\leq s\leq v_-^{-1}(n)+v_+^{-1}(n)+1}t^{2s}(\kappa_u(-v_-^{-1}(n),v_+^{-1}(n)))^s}+\frac{C}{\varphi_u^{-1}(1/t)}\nonumber\\
&\leq& \frac{C}{n}\sqrt{1+\sum_{1\leq s\leq v_-^{-1}(n)+v_+^{-1}(n)+1}t^{2s}\varphi_u^{2s}(n)}+\frac{C}{\varphi_u^{-1}(1/t)}.\nonumber\end{eqnarray}

\medskip
\noindent
Choose $n=\varphi_u^{-1}(1/2t)$. In particular, $\varphi_u(n)\leq1/(2t)$. This gives~:

$$|1-\chi_D(ut)|\leq \frac{C}{\varphi_u^{-1}(1/(2t))}\sqrt{1+\sum_{s\geq1}(1/2)^{2s}}+\frac{C}{\varphi_u^{-1}(1/t)}.$$

\medskip
\noindent
By lemma \ref{young}, there is a constant $C'$ independent on $u\in S_+^{d-1}$ so that for small $t>0$~:

$$|1-\chi_D(ut)|\leq \frac{C'}{\varphi_u^{-1}(1/t)}.$$

\medskip
\noindent
This concludes the proof of the proposition. \fin

\subsection{Conclusion}

\noindent
- {\it Theorem \ref{princip}, corollary \ref{principcor} and proposition \ref{antiscor}}. 
By propositions \ref{re}, \ref{mod} and theorem \ref{spitz}, using that $\mbox{Re}(1/a)=\mbox{Re}(a)/|a|^2$, the random walk is recurrent if and only if, for some $\eta>0$~:

\begin{equation}
\label{stereo}
\int_{(u,t)\in S_+^{d-1}\times(0,\eta)}\frac{(\varphi_u^{-1}(1/t))^2}{\varphi_{u,+}^{-1}(1/t)}~t^{d-1}dudt=+\infty.
\end{equation}

\medskip
\noindent
For fixed $u\in S_+^{d-1}$, we cut the interval $(0,\eta)$ in the contiguous intervals $[1/(n+1),1/n]$, $n\geq n_0$. The latter have length of order $1/n^2$, so using finally lemma \ref{young}, the condition is equivalent to the one given in the statement of theorem \ref{princip}.

\medskip
Concerning proposition \ref{antiscor}, we first show in the antisymmetric case that $\varphi_u^{-1}$ and $\varphi_{u,++}^{-1}$ have the same size, uniformly in $u\in S_+^{d-1}$. By lemma \ref{young}, it is enough to show that $\varphi_u\leq C\varphi_{u,++}$. Observe that $p_0=q_0$ and~:

\begin{eqnarray}
\varphi_{u}^2(n)&=&\varphi_{u,+}^2(n)+\sum_{-v_-^{-1}(n)\leq k\leq 0\leq l\leq v_+^{-1}(n)}T_k^l(u)\nonumber\\
&=&\varphi_{u,+}^2(n)+\sum_{0\leq k,l\leq v_+^{-1}(n)}T_{\min(k,l)+1}^{\max(k,l)}(u)\leq 4\varphi_{u,+}^2(n)\leq 8\varphi_{u,++}^2(n).\nonumber
\end{eqnarray}

\smallskip
This completes the proof of this claim.

\medskip
Concerning Corollary \ref{principcor}, we always have $\varphi_u^{-1}\leq \varphi_{u,+}^{-1}$. Then~:

\begin{equation}
\label{stereo2}
\int_{S_+^{d-1}\times(0,\eta)}\frac{(\varphi_u^{-1}(1/t))^2}{\varphi_{u,+}^{-1}(1/t)}~t^{d-1}dudt\leq\int_{S_+^{d-1}\times(0,\eta)}\varphi_u^{-1}(1/t)t^{d-1}dudt.
\end{equation}

\noindent
\smallskip
In the antisymmetric case, both integrals have the same order. To complete the proofs of proposition \ref{antiscor} and corollary \ref{principcor}, we just need to show that the second term has the right order. For fixed $u\in S_+^{d-1}$, up to decreasing $\eta>0$, also taking $n_0$ independent on $u\in S_+^{d-1}$ (as $0<\alpha\leq \varphi_u(1)\leq\beta $, for constants $\alpha$ and $\beta$, independent on $u\in S_+^{d-1}$)~:

$$\sum_{n\geq n_0}\int_{1/\varphi_u(n+1)}^{1/\varphi_u(n)}\varphi_u^{-1}(1/t)t^{d-1}~dt\leq \int_0^{\eta}\varphi_u^{-1}(1/t)t^{d-1}~dt\leq\sum_{n\geq 1}\int_{1/\varphi_u(n+1)}^{1/\varphi_u(n)}\varphi_u^{-1}(1/t)t^{d-1}~dt.$$

\noindent
\smallskip
On each domain $(1/\varphi_u(n+1),1/\varphi_u(n))$, we have $\varphi_u^{-1}(1/t)=n$. Hence $\int_0^{\eta}\varphi_u^{-1}(1/t)t^{d-1}~dt$ has exact order~:

\begin{eqnarray}
\sum_{n\geq 1}n\int_{1/\varphi_u(n+1)}^{1/\varphi_u(n)}t^{d-1}~dt&=&\frac{1}{d}\sum_{n\geq 1}n\({\frac{1}{(\varphi_u(n))^d}-\frac{1}{(\varphi_u(n+1})^d}\)\nonumber\\
&=&\frac{1}{d}\lim_{N\rightarrow+\infty}\sum_{n=2}^N\frac{1}{(\varphi_u(n))^d}+\({\frac{1}{(\varphi_u(1)^d}}-\frac{N}{(\varphi_u(N+1))^d}\).\nonumber\\
&=&\frac{1}{d}\lim_{N\rightarrow+\infty}\sum_{n=1}^N\({\frac{1}{(\varphi_u(n))^d}}-\frac{1}{(\varphi_u(N+1))^d}\).\nonumber
\end{eqnarray}

\smallskip
\noindent
Remark that the right-hand side is bounded from above by $(1/d)\sum_{n\geq1}1/(\varphi_u(n))^d$. Hence~:

- if $\int_{u\in S_+^{d-1}}\sum_{n\geq1}(1/\varphi_u(n))^d<+\infty$, then the left-hand side in \eqref{stereo2} is finite. 

\smallskip
- if $\int_{u\in S_+^{d-1}}\sum_{n\geq1}(1/\varphi_u(n))^d=+\infty$, using at the end Fatou's lemma~:

\begin{eqnarray}
\int_{S_+^{d-1}\times(0,\eta)}\varphi_u^{-1}(1/t)t^{d-1}dudt&\geq& C\int_{S_+^{d-1}}\lim_{N\rightarrow+\infty}\sum_{n=1}^N\({\frac{1}{(\varphi_u(n))^d}}-\frac{1}{(\varphi_u(N+1))^d}\)~du\nonumber\\
&\geq &C\int_{S_+^{d-1}}\sum_{n\geq1}\frac{1}{(\varphi_u(n))^d}~du=+\infty.\nonumber\end{eqnarray}

\noindent
This completes the proofs of corollary \ref{principcor} and of the first part of proposition \ref{antiscor}. To complete the proof of the latter, we take $m_0=0$ and $m_n=-m_{-n}=c\not=0$, $n\geq1$. Then $\varphi_{u,++}^2(n)$ has immediately the same order as~:

$$nw_+\circ v_+^{-1}(n)+(c.u)^2\sum_{1\leq k<l\leq v_+^{-1}(n)}\rho_k\rho_l\({\sum_{k\leq s\leq l}1/\rho_s}\)^2.$$

\noindent
We still denote by $\varphi_{u,++}^2(n)$ this quantity. Suppose now that $c_1n^{\alpha}\leq \rho_n\leq c_2n^{\alpha}$, $n\geq0$. We reason up to multiplicative constants, using the notation $\asymp$.

\medskip
-- If $\alpha<-1$, the random walk is transient, as $(v_+(n))$ is bounded.

\medskip
-- If $\alpha=-1$, then $w_+(n)\asymp n^2$ and $v_+(n)\asymp \ln n$. As a result, $\varphi_{u,++}^2(n)\geq ne^{cn}$, for some $c>0$, giving transience.

\medskip
-- Suppose that $-1<\alpha<1$. We show transience.  We have $w_+(n)\asymp n^{1-\alpha}$ and $v_+(n)\asymp n^{1+\alpha}$.

\begin{eqnarray}
\sum_{1\leq k<l\leq n}\rho_k\rho_l\({\sum_{k\leq s\leq l}1/\rho_s}\)^2&\asymp& \int_{1\leq x\leq y\leq n}x^{\alpha}y^{\alpha}\({\int_x^yt^{-\alpha}~dt}\)^2dxdy\nonumber\\
&\asymp& n^{2\alpha+2}\int_{1/n\leq x\leq y\leq 1}x^{\alpha}y^{\alpha}\({\int_{nx}^{ny}t^{-\alpha}~dt}\)^2dxdy\nonumber\\
&\asymp& n^{4}\int_{1/n\leq x\leq y\leq 1}x^{\alpha}y^{\alpha}\({\int_{x}^{y}t^{-\alpha}~dt}\)^2dxdy\asymp n^4.\nonumber\end{eqnarray}

\noindent
As a result $\varphi_{u,++}^2(n)\asymp n^{1+(1-\alpha)/(1+\alpha)}+(c.u)n^{4/(1+\alpha)}$, so $\varphi_{u,++}(n)\asymp n^{1/(1+\alpha)}+(c.u)n^{2/(1+\alpha)}$. We obtain that when $d=1$, $\varphi_{u,++}(n)\asymp n^{2/(1+\alpha)}$ and when $d\geq2$, $\varphi^d_{u,++}(n)\geq cn^{d/(1+\alpha)}$. As the exponents are $>1$ in each case, the random walk is transient, from corollary \ref{principcor}.

\medskip
-- Suppose next that $\alpha>1$. Then $w_+(n)\asymp 1$, $v_+(n)\asymp n^{1+\alpha}$. If $d=1$, then $\varphi_{u,++}^2(n)\leq C(n+(c,u)^2n^2)$, so $\varphi_{u,++}(n)=O(n)$ and the random walk is recurrent. When $d=2$~:

$$\varphi_{u,++}^2(n)\asymp n+(c.u)^2 \int_{1\leq x\leq y\leq n}x^{\alpha}y^{\alpha}\({\int_x^yt^{-\alpha}~dt}\)^2dxdy.$$

\noindent
The second term can be written as~:

\begin{eqnarray}
&&\int_1^nx^{\alpha}~dx\int_1^nx^{\alpha}\({\int_x^{+\infty}t^{-\alpha}~dt}\)^2dx-\({\int_1^nx^{\alpha}\int_x^{+\infty}t^{-\alpha}dt~dx}\)^2\nonumber\\
&\asymp&\({\int_1^nx^{\alpha}~dx}\)\({\int_1^nx^{2-\alpha}}\)-\({\int_1^nx~dx}\)^2.\nonumber
\end{eqnarray}

\noindent
Let $1<\alpha<3$. Then this term is equivalent to $((\alpha-1)^2/(\alpha+1)(3-\alpha))n^4$. As a result~:

$$\varphi_{u,++}^2(n)\asymp n+(c.u)^2n^{4/(1+\alpha)}.$$

\medskip
\noindent
In order to show transience we need to control the following quantity~:

$$\sum_{n\geq1}\int_{u\in S^1_+}\frac{1}{\varphi^2_{u,++}(n)}\asymp\sum_{n\geq1}\int_0^{\pi/2}\frac{1}{n+\theta^2n^{4/(1+\alpha)}}~d\theta.$$

\noindent
Setting $\theta=n^{1/2-2/(1+\alpha)}x$, it remains~:

$$\sum_{n\geq1}\frac{1}{n}n^{1/2-2/(1+\alpha)}\int_0^{(\pi/2)n^{2/(1+\alpha)-1/2}}\frac{1}{1+x^2}~dx\asymp\sum_{n\geq1}\frac{1}{n^{1/2+2/(1+\alpha)}}<+\infty,$$

\noindent
as $1/2+2/(1+\alpha)>1$. If $\alpha=3$, then $\varphi_{u,++}^2(n)\asymp n+(c.u)^2n\ln n\leq Cn\ln n$. When $\alpha>3$,  $\varphi_{u,++}^2(n)\asymp n+(c.u)^2n\leq Cn$. In any case $\sum_{1\geq1}(1/\varphi_{u,++}^2(n))=+\infty$, giving recurrence.

\medskip
-- If $\alpha=1$, then $w_+(n)\asymp \ln n$, $v_+(n)\asymp n^2$. When $d=1$, $\varphi_{u,++}^2(n)\leq C(n\ln n+n^2(\ln n)^2)$, so $\varphi_{u,++}(n)=O(n\ln n)$ and the random walk is recurrent. When $d=2$, notice that~:
 
$$\varphi_{u,++}^2(n)\geq K(n\ln n+(c.u)^2n^2).$$

\noindent
In order to show transience, we just need to prove the finiteness of~:

$$\sum_{n\geq1}\int_0^{\pi/2}\frac{1}{n\ln n+\theta^2n^2}~d\theta=\sum_{n\geq 1}\int_0^{(\pi/2)\sqrt{n/\ln n}}\frac{dx}{1+x^2}\sqrt{(\ln n)/n}\frac{1}{n\ln n}<+\infty.$$

\noindent
This completes the proof of the proposition.

\medskip
\noindent
- {\it Proposition \ref{halfpipe}}. Let $\tilde{\varphi}$ and $\tilde{\varphi_+}$ be the functions corresponding to the case when $m_n=1$, $n\in\Z$. Set $D=\sum_{n\in\Z}(m_n/\rho_n)$. Observe that one always has in the present situation~:

$$\varphi^2(n)=\sum_{-v_-^{-1}(n)\leq k\leq l\leq v_+^{-1}(n)}\rho_{k-1}\rho_l\({\sum_{k\leq r\leq l}(m_r/\rho_r)}\)^2\leq C\sum_{-v_-^{-1}(n)\leq k\leq l\leq v_+^{-1}(n)}\rho_{k-1}\rho_l\leq Cn^2.$$

\noindent
Let $A>0$ so that $v_{\pm}(n+1)\leq (A/2) v_{\pm}(n)$. In the case when $D\not=0$, one also has~:

\begin{eqnarray}
\varphi^2(n)&\geq& \sum_{-v_-^{-1}(n)\leq k\leq -v_-^{-1}(n/A),v_+^{-1}(n/A)\leq l\leq v_+^{-1}(n)}\rho_{k-1}\rho_l\({\sum_{k\leq r\leq l}(m_r/\rho_r)}\)^2\nonumber\\
&\geq&(n/A)(n/A)(D/2)^2,\nonumber
\end{eqnarray}

\noindent
if $n$ is large enough. As a result $\varphi(n)$ and $\varphi^{-1}(n)$ have order $n$. The same is true for $\tilde{\varphi}$. It remains to show the finiteness of~:

$$\sum_{n\geq1}\frac{1}{\varphi_+^{-1}(n)}\leq C\sum_{n\geq1}\frac{1}{\tilde{\varphi_+}^{-1}(n)}\leq C\sum_{n\geq1}\frac{(\tilde{\varphi}^{-1}(n))^2}{n^2\tilde{\varphi_+}^{-1}(n)}<+\infty,$$

\noindent
because the random walk is obviously transient when $m_n=1$, $n\in\Z$.

\medskip
When $D=0$ and in the antisymmetric case, by corollary \ref{principcor} the criterion reduces to~:

$$\sum_{n\geq1}\frac{1}{\varphi_+(n)}\geq\sum_{n\geq1}\frac{1}{\varphi(n)}\geq c\sum_{n\geq1}\frac{1}{n}=+\infty.$$

\noindent
The random walk is recurrent. This completes the proof of the proposition.
\subsection{Remarks}

It seems necessary to interpret the recurrence criterion in order to use it in practice. When $\mu_n=\delta_c$, with $c\not=0$, the integral in the criterion is finite, because the random walk is trivially transient. How does one may see it directly ? The question is not clear, even for $d=1$ and $c=1$.

\medskip
It would be interesting to consider the case when the $(p_n,q_n,m_n)$ are a typical realization of an i.i.d. process with $m_n$ independent of $(p_n,q_n)$, $\E(\log (p_n/q_n))=0$, $\mbox{Var}(\log (p_n/q_n))>0$, $\E(m_n)=0$ and $\mbox{var}(m_n)>0$. One needs first of all to study in detail $(v_{\pm}(n))$. The random walk is without any doubt transient.

\medskip
It would also be of interest to consider the analogous model in $\Z\times\Z^2$ in a $\Z^2$-invariant environment. If following the main strategy, the main difficulty in proving a characterization of the asymptotical behaviour is to detail the distribution of the local time during an excursion of simple random walk in the plane. There is no tree-structure behind, but a complicated graph with loops. A first step in this direction seems to be the following model in the plane~:

$$P_{(m,n),(m,n\pm1)}=1/4,~\P_{(m,n),(m+1,n)}=p(m,n)/2,~\P_{(m,n),(m-1,n)}=q(m,n)/2,$$

\smallskip
\noindent
with $p(m,n)+q(m,n)=1$, for example making some hypothesis of stochastic homogeneity on the $(q(m,n),p(m,n))_{(m,n)\in\Z^2}$. The vertical component being recurrent, one may study the subsequence of return times on the horizontal axis. This random walk is a one-dimensional random walk in random medium with unbounded jumps. Very few results are known on such a random walk, cf Andjel \cite{andjel}, and they suppose the jump integrable, which is not the case here. The very first step in the proof (lemma \ref{petri}) is already not clear.


\providecommand{\bysame}{\leavevmode\hbox to3em{\hrulefill}\thinspace}

\bigskip
{\small{\sc{Laboratoire d'Analyse et de Math\'ematiques Appliqu\'ees,
  Universit\'e Paris-Est, Facult\'e des Sciences et Technologies, 61, avenue du G\'en\'eral de Gaulle, 94010 Cr\'eteil Cedex, FRANCE}}}

\it{E-mail address~:} {\sf julien.bremont@u-pec.fr}

\end{document}